\documentclass[proc]{edpsmath}
\usepackage[utf8]{inputenc}
\usepackage[T1]{fontenc}

\usepackage{amsmath,amssymb,amsthm}
\usepackage{graphicx}
\usepackage{hyperref}
\usepackage{physics}
\usepackage{mathtools}
\usepackage{float}

\newcommand{\bR}{\mathbb{R}}
\newcommand{\bb}{\mathbf{b}}

\newcommand{\be}{\mathbf{e}}
\newcommand{\vp}{v_\parallel}
\newcommand{\bx}{\mathbf{x}}

\newcommand{\bu}{\mathbf{u}}
\newcommand{\bB}{\mathbf{B}}
\newcommand{\Dt}{\Delta t \;}

\newcommand{\dbx}{\dot{\mathbf{x}}}
\newcommand{\diff}[2]{\frac{\text{d} #1}{\text{d} #2}}
\newcommand{\pa}[2]{\frac{\partial #1}{\partial #2}}
\newcommand{\ap}{a_\parallel}
\newcommand{\mean}[1]{\langle #1 \rangle}
\newcommand{\Bap}{B_\parallel^\ast}
\newcommand{\bracket}[2]{\left\{ #1, #2 \right\}}
\newcommand{\gcbracket}[2]{\left\{ #1, #2 \right\}_{\text{g.c.}}}
\newcommand{\dt}{\frac{\text{d}}{\text{d} t}}

\renewcommand{\d}{\;\text{d}}

\usepackage{tikz}
\usetikzlibrary{cd}
\tikzcdset{arrow style=tikz
	, diagrams={>={Stealth[length=1.2mm, width=1mm]}}
}
\usetikzlibrary{positioning}
\usetikzlibrary{arrows.meta}
\usetikzlibrary{decorations.pathreplacing}
\usetikzlibrary{shapes.symbols}
\usetikzlibrary{shapes.geometric}

\begin{document}

\title{Splitting Scheme for Gyro-kinetic Equations with Semi-Lagrangian and Arakawa Substeps} 

\author{Dominik Bell}\address{Max-Planck-Institut für Plasmaphysik, Garching, Germany; \email{dominik.bell@ipp.mpg.de\ \&\ martin.campos-pinto@ipp.mpg.de\ \&\ frederik.schnack@ipp.mpg.de\ \&\ sonnen@ipp.mpg.de}}\secondaddress{Technische Universität München, Zentrum Mathematik, Garching, Germany}
\author{Martin Campos Pinto}\sameaddress{1}
\author{Davor Kumozec}\address{Faculty of Sciences, University of Novi Sad, Serbia; \email{davor.kumozec@dmi.uns.ac.rs}}
\author{Frederik Schnack}\sameaddress{1}
\author{Emily Bourne}\address{CEA\unskip, IRFM\unskip, Saint-Paul-les-Durance\unskip, F-13108\unskip, France; \email{emily.bourne@epfl.ch}}
\author{Eric Sonnendrücker}\sameaddress{1,2}

\begin{abstract}
	The gyro-kinetic model is an approximation of the Vlasov-Maxwell system in a strongly magnetized magnetic field. We propose a new algorithm for solving it combining the Semi-Lagrangian (SL) method and the Arakawa (AKW) scheme with a time-integrator. Both methods are successfully used in practice for different kinds of applications, in our case, we combine them by first decomposing the problem into a fast (parallel) and a slow (perpendicular) dynamical system. The SL approach and the AKW scheme will be used to solve respectively the fast and the slow subsystems. Compared to the scheme in \cite{pygyro_code}, where the entire model is solved using only the SL method, our goal is to replace the method used in the slow subsystem by the AKW scheme, in order to improve the conservation of the physical constants.
\end{abstract}

\begin{resume}
	Le modèle gyro-cinétique est une approximation du système de Vlasov-Maxwell dans un champ magnétique fortement magnétisé. Nous proposons un nouvel algorithme pour le résoudre en combinant la méthode Semi-Lagrangienne (SL) et le schéma d'Arakawa (AKW) avec un intégrateur temporel. Les deux méthodes sont utilisées avec succès dans la pratique pour différents types d'applications. Dans notre cas, nous les combinons en décomposant d'abord le problème en un système dynamique rapide (parallèle) et un système dynamique lent (perpendiculaire). L'approche SL et le schéma AKW seront utilisés pour résoudre respectivement les sous-systèmes rapide et lent. Par rapport au schéma de \cite{pygyro_code}, où le modèle entier est résolu en utilisant uniquement la méthode SL, notre objectif est de remplacer la méthode utilisée dans le sous-système lent par le schéma AKW, afin d'améliorer la conservation des constantes physiques.
\end{resume}

\maketitle
\newpage
\tableofcontents


\section{Introduction}

Gyro-kinetic models have become a very popular choice for simulating plasmas since they reduce the phase space compared to the Vlasov equation and are thus more computationally viable while still describing the important physics. This reduction is done by averaging out the fast motion of particles around the magnetic field lines, called gyration.

In the works \cite{idomura2007} and \cite{idomura2008conservative} the authors use a Morinishi finite differences (FD) scheme \cite{Morinishi1998FullyCH} for the advection equations, which conserves the mass and $L^2$-norm of the distribution function. In \cite{crouseilles2018exponential}, the gyro-kinetic equation is split; the linear part of the advection is solved using Fourier techniques while the non-linear part is discretized using the FD Arakawa scheme from \cite{Arakawa_1966} of order 2, which aims at preserving the kinetic energy and the square vorticity. The transport in time then uses an exponential integrator. In \cite{Latu_2017} the model is split into a fast and a slow subsystem for all of which a backward Semi-Lagrangian scheme is used.

This paper follows the latter method but replaces the Semi-Lagrangian scheme in the slow subsystem (which does not involve the magnetic field) by an Arakawa scheme of order 4. The motivation behind this is to improve conservation properties for the most turbulent substep at the expense of either using a costly implicit time integrator, or an explicit one which is constrained by a CFL condition. For this method we test different orders and boundary conditions. This spatial discretization will be combined with either an implicit Crank-Nicolson integrator or an explicit Runge-Kutta scheme of order 4. The implementation proceeds on top of the \texttt{PyGyro} code \cite{pygyro_code} which is a \textit{Python} implementation of the gyro-kinetic model presented in \cite{Latu_2017} with the aim to replicate test cases of the \texttt{GYSELA} code from \cite{GRANDGIRARD2006395} and more recently \cite{Grandgirard_CPC2016}.

After an introduction of the gyro-kinetic equation in Section \ref{sec:gk-model}, we will describe the splitting ansatz of dividing the model in a slow and a fast subsystem, as well as the two aforementioned schemes in Section \ref{sec:splitting_discretization}. In section \ref{sec:num_exp}, numerical experiments will be presented, first as a verification of the Arakawa scheme in a 2-dimensional test case, and after that results for the full gyro-kinetic model. Concluding with Section \ref{sec:conclusion}, we discuss the benefits of this approach and perspectives.


\section{The Gyro-kinetic Model}
\label{sec:gk-model}

In a Tokamak, the dynamics of particles consists of a slow motion along the magnetic field lines superimposed with a fast gyration around the magnetic field lines. This fast motion can be averaged out to reduce the dimension of phase space in order to make the models computationally more tractable, while keeping most of the important physics. The resulting models are called gyro-kinetic and will be described in the following.\\

The model we will look at is defined by the Lie-transformed, low-frequency particle Lagrangian $L$, where we follow the derivation in \cite{Bottino_Sonnendrucker_2015}, \cite{emily} and \cite{Latu_2017}. Given a static magnetic field $\bB$, its intensity $B = \norm{\bB}$ and its direction $\bb = \bB / B$, the particle charge $q \not= 0$ and the particle mass $m >0$,  
we are able to write the Lagrangian $L$ as
\begin{equation}\label{Lagrangian}
	L(t, \bx, \vp, \mu, \dot{\bx}, \dot{v}_\parallel, \dot{\mu}) = \left(\nabla \times q\bB + m \vp \bb \right)\cdot\dbx + \frac{m}{q} \mu \dot{\theta} - H(t, \bx, \vp),
\end{equation}
where $\bx \in \Omega \subseteq \bR^3$ is the position of the gyro-centre, $\vp \in \bR$ is the velocity parallel to the magnetic field lines, $\mu$ is the modified magnetic moment and $\theta$ the angle of cylindrical coordinates. The Hamiltonian $H$ will be introduced shortly.
Looking at the equation of motion of $\theta$, i.e. its Euler-Lagrange equation
\begin{equation}
	0 = \frac{\dd}{\dd t} \left( \frac{\partial L}{\partial \dot \theta} \right) - \left(
		\frac{\partial L}{\partial \theta}\right) =  \frac{\dd}{\dd t} \frac{m}{q} \mu,
\end{equation}
we can immediately conclude that $\mu$ is an exact invariant of the system, i.e.
\begin{equation}
	\diff{}{t} \mu = 0,
\end{equation}
and thus the phase-space is only $4$-dimensional.
The gyro-kinetic equation describing the gyro-centre distribution function $f=f(t,\bx, \vp)$ is of the form
\begin{equation}\label{drift-kinetic model}
	\pa{f}{t} + \bu \cdot \nabla f + \ap \pa{f}{\vp} = 0,
\end{equation}
which describes the positions of a collection of identical particles and whose exact solution is constant along the trajectories $(\bx(t), \vp(t))$ in the phase-space, i.e. 
\begin{equation}
	\frac{\dd}{\dd t} f(t,\bx(t), \vp(t)) = 0.
\end{equation}
In order to determine the equations of motion for ${\dd \bx}/{\dd t} = \bu$ and ${\dd \vp}/{\dd t} = \ap$, we look at the remaining Euler-Lagrange equations. Simplifying notations by defining
\begin{equation}
	\bB^\ast  \coloneqq \bB + \frac{m}{q}\vp \nabla \times \bb , \qquad
	\Bap  \coloneqq \bb \cdot \bB^\ast = B + \frac{m \vp}{q B} \bb \cdot \left( \nabla \times \bB \right),
\end{equation}
one can derive the characteristic trajectories (see \cite{Bottino_Sonnendrucker_2015} for a detailed derivation) from the remaining Euler-Lagrange equations of \eqref{Lagrangian}, which yield
\begin{subequations}
	\begin{align}
		\bu  & = \frac{1}{\Bap} \left( \frac{1}{m} \pa{H}{\vp} \bB^\ast + \frac{1}{q} \bb \times \nabla H \right) \label{eom for bu}, \\
		\ap & = \frac{1}{\Bap} \left( -\frac{1}{m} \bB^\ast \cdot  \nabla H  \right). \label{eom for ap}
	\end{align}
\end{subequations}
As noted in \cite{Latu_2017}, the phase space is divergence-free, i.e.
\begin{equation}
	\nabla \cdot \bu + \pa{\ap}{\vp} = 0,
\end{equation}
thus we can rewrite \eqref{drift-kinetic model} in conservative form
\begin{equation}\label{conservation1}
	\pa{}{t}\left(\Bap f\right) + \nabla\cdot\left( \Bap \bu f \right) + \pa{}{\vp} \left( \Bap \ap f \right) = 0.
\end{equation}

The equations of motion form a Hamiltonian system (see \cite{idomura2008conservative} for details) with the electrostatic gyro-centre Hamiltonian in the zero-Larmor-radius limit
\begin{equation}
	H(t, \bx, \vp) = \frac{1}{2} m \vp^2 + \mu B(\bx) + q \phi (t,\bx).
\end{equation}
Furthermore, we need a bracket for the Hamiltonian system which is the guiding-centre Poisson bracket
\begin{subequations}
	\begin{align}
		\gcbracket{F}{G} & = \frac{\bB}{m \Bap} \cdot \left( \left(\nabla F\right) \pa{G}{\vp} - \left(\nabla G\right) \pa{F}{\vp} \right) \label{fast-poisson} \\
		& \qquad + \frac{\vp}{q \Bap} \left( \nabla \times \bb \right) \cdot \left( \left(\nabla F\right) \pa{G}{\vp} - \left(\nabla G\right) \pa{G}{\vp} \right) \\
		& \qquad - \frac{1}{q \Bap} \bb \cdot \left[\left(\nabla F\right) \times \left(\nabla G\right)\right].
	\end{align}
\end{subequations}
We now split the Poisson bracket in parts containing $\bB$ (i.e. \eqref{fast-poisson}) which we expect to have fast dynamics, and other terms which will be called the slow subsystem.
Thus the equations of motion for the fast and slow subsystems, with variables $\bu = \bu_f + \bu_s$ and $a_{\parallel} = a_{\parallel, f} + a_{\parallel, s}$, read:
\begin{subequations}
	\begin{align}
		(\text{fast}) = & \left\{ \begin{aligned}
			\bu_f & = \frac{1}{\Bap} \frac{1}{m} \pa{H}{\vp} \bB, \\
			a_{\parallel, f} & = -\frac{1}{\Bap} \frac{1}{m} \bB \cdot \left( \nabla H \right),
		\end{aligned} \right. \label{split step fast}\\
		(\text{slow}) = &  \left\{ \begin{aligned}
			\bu_s & = \frac{1}{\Bap} \frac{1}{q} \left( \pa{H}{\vp} \vp \nabla \times \bb + \bb \times \left(\nabla H\right) \right), \\
			a_{\parallel, s} & = - \frac{a}{\Bap} \frac{1}{q} \vp \left( \nabla \times \bb \right) \cdot \left( \nabla H \right).
		\end{aligned} \right. \label{split step slow}
	\end{align}
\end{subequations}
This splitting is of particular interest when introducing a phase-space discretization, where the fast system trajectories may travel across many cells, which introduces a higher CFL number and thus needs a time-integration that handles this well, while the treatment of the slow subsystem naturally is less restrictive.

For a constant and uniform background magnetic field $\bB = B \hat \be_z$, where $\hat \be_z$ is the unit-vector pointing in the $z$-variable direction, the model simplifies to drift-kinetic equations and the subsystems become
\begin{subequations}
	\begin{align}
		(\text{fast}) = & \left\{ \begin{aligned}
			\bu_f & = {\vp} \bb, \\
			a_{\parallel, f} & = -\frac{q}{m} \bb \cdot \nabla \phi,
		\end{aligned} \right. \label{split step fast const B} \\
		(\text{slow}) = &  \left\{ \begin{aligned}
			\bu_s & = \frac{\bb \times \nabla \phi}{B}, \\
			a_{\parallel, s} & = 0.
		\end{aligned} \right. \label{split step slow const B}
	\end{align}
\end{subequations}
Continuing this simplification, we are able to rewrite the space variables to cylindrical coordinates as are introduced in \cite{Latu_2017}. Thus, we look for the distribution function $f = f(t, r, \theta, z, v_\parallel)$ satisfying
\begin{equation}
	\partial_t f + \{\phi, f\} + v_\parallel \nabla_\parallel f - \nabla_\parallel \phi \partial_{v_\parallel}f = 0, \label{eq:GK_model}
\end{equation}
where $\nabla_\parallel = \bb \cdot \nabla$ and the bracket is transformed to polar coordinates, which reads, given the toroidal magnetic field $B_0$,
\begin{equation}
	\{\phi, f \} = \frac{1}{rB_0}\partial_r \phi\partial_\theta f -\frac{1}{rB_0}\partial_\theta \phi\partial_r f. \label{eq:bracket_in_polar}
\end{equation}
Since the plasma is quasi-neutral this equation is complemented by solving an elliptic quasi-neutrality (QN) equation for the self-consistent potential $\phi = \phi(t, r, \theta, z)$, i.e. solving for a given temperature profile $T_e$
\begin{equation}
	- \left[\partial_r^2 \phi + \left( \frac{1}{r} + \frac{\partial_r n_0}{n_0}\right)\partial_r \phi + \frac{1}{r^2} \partial_\theta^2 \phi \right] + \frac{1}{T_e} \phi = \int_{-\infty}^{\infty} (f - f_\text{eq}) \ \dd v_\parallel, \label{eq:qn}
\end{equation}
where the given (radial symmetric) equilibrium function $f_\text{eq}$ is a Gaussian and $n_0$ is a radial profile, which is the integral over $\vp$ of the equilibrium function. The initial distribution function $f(t=0, r, \theta, z, v_\parallel)$ is defined to be the equilibrium with a perturbation in some modes, which will result in an observable turbulence behaviour of the system after a certain amount of time. For more details and the exact definitions and constants, we refer to Section 4 in \cite{Latu_2017}. Normally these equations are defined on an infinitely long cylinder of some radius, but in order to discretize them, we cut the domain down to a finite cylinder with a hole in the middle and reduce the velocity space, such that   $(r, \theta, z, \vp) \in [r_\text{min}, r_\text{max}]\times[0, 2\pi]\times[0, 2\pi R_0] \times [-v_\text{max}, v_\text{max}]$, with parameters that will be made more precise in Section \ref{sec:num_exp}. To complete this system of equations, we briefly discuss the boundary conditions of this reduction. The distribution function $f$ is periodic along $\theta$, $z$, while having homogeneous Dirichlet boundaries in the $\vp$-direction. In the radial direction $r$, we assume the values are given by an outside equilibrium function. For the potential $\phi$, we assume periodic boundary conditions along $\theta$ and $z$. In $r$-direction, we decompose $\phi$ into Fourier modes at $r_\text{min}$, taking homogeneous Neumann boundary conditions for the zeroth mode and homogeneous Dirichlet boundary conditions for the others. At $r_\text{max}$, we just take plain homogeneous Dirichlet boundary conditions. \\

Lastly, we want to take a look at physical properties of these equations. Since we have a transport equation in conservative form \eqref{conservation1} that conserves arbitrary functions of $f$  along non-linear characteristic trajectories, we can write the Casimir equation
\begin{equation}
	\dt C(f) + \bracket{C(f)}{H} =0,
\end{equation}
where the Casimir invariant $\int C( f ) $ is conserved. Therefore, the system has an infinite number of conserved quantities such as the particle number, the $L^1$-norm $\int f $, or the $L^2$-norm $\int f^2 $. In addition, the gyro-kinetic equations conserve the total energy $H = E_k + E_p$:
\begin{subequations}
	\begin{align}
		E_k &= \int f(\bx, \vp) \left[\frac{1}{2} m \vp^2 + \mu B(\bx)\right] \d \bx \d \vp \, , \\
		E_p &= \int q \mean{\phi}_\alpha (t,\bx) f(\bx, \vp) \d \bx \d \vp \, ,
	\end{align}
\end{subequations}
where $E_k$ is the kinetic energy and $E_f$ is the potential energy. For more details we refer the reader to \cite{idomura2008conservative}.


\section{Discretization via Operator Splitting}
\label{sec:splitting_discretization}

Operator splitting has proven to be an effective method for the time-integration of ordinary differential equations (ODEs), whose vector-field can be written as a sum of simpler terms. This is of special interest for geometric integration, that is, if the underlying problem has geometric properties that should be conserved by the integrator and can be enforced more easily for each split step. Usually, the combined computational cost for the integration of the split steps is lower compared to integrating the full ODE, albeit at the expense of introducing an additional approximation error. A broad overview of these methods is given in \cite{mclachlan_quispel_2002}, also containing the second-order Strang and first-order Lie splitting, that we will use in the following.\\

In the SL scheme from \cite{pygyro_code}, the (Hamiltonian) operator splitting is applied on the full screw-pinch model equation, that has a separable Hamiltonian
\begin{equation}
	H = \frac{1}{2} m \vp^2 + q \phi
\end{equation}
and describes the time-evolution of the distribution function $f$, i.e.
\begin{equation}
 \partial_t f + \{\phi, f \} + v_\parallel \nabla_\parallel f - \nabla_\parallel \phi\,\, \partial_{v_\parallel} f = 0,
\end{equation}
We can split this by first applying the above discussed splitting of the Poisson bracket into fast and slow parts as described above, and then splitting the Hamiltonian into its two parts in the fast subsystem, yielding
\begin{subequations}
	\begin{align}
		\partial_t f + v_\parallel \nabla_\parallel f & = 0, && \text{(Advection on flux surface),} & \label{eq:adv_flux} \\
		\partial_t f + \nabla_\parallel \phi\,\, \partial_{v_{\parallel}} f & = 0, && (v\text{-parallel advection),} & \label{eq:adv_par} \\
		\partial_t f + \{\phi, f\} & = 0, && \text{(Advection on poloidal plane),} \label{eq:adv_poloidal} &
	\end{align}
\end{subequations}
where $\{\phi,f\}$ is the Poisson bracket in polar coordinates defined in \eqref{eq:bracket_in_polar}. Splitting the Hamiltonian immediately and using the whole gyro-centre bracket would yield only 2 equations, namely \eqref{eq:adv_par} and \eqref{eq:adv_poloidal} would be kept together.\\

As described in \cite{emily} and \cite{Latu_2017}, we then obtain solution operators to each sub-system by the SL method, which will shortly be introduced in the next chapter. For now, we denote them by $A,B$ and $C$ solving the equations \eqref{eq:adv_flux}, \eqref{eq:adv_par} and \eqref{eq:adv_poloidal} for a given time-step and potential $\phi$ respectively.
For a time step $\Delta \tau \in \bR$, the first-order Lie splitting of $f^n = f(t = n \Delta \tau)$ reads
\begin{equation}
	f^{n+1/2} = C\left(\frac{\Delta \tau}{2}\right) B\left(\frac{\Delta \tau}{2}\right) A\left(\frac{\Delta \tau}{2}\right) f^n, \label{eq:Lie}
\end{equation}
while the second-order Strang splitting is given by
\begin{equation}
	f^{n+1} = A\left(\frac{\Delta \tau}{2}\right) B\left(\frac{\Delta \tau}{2}\right) C\left(\Delta \tau\right) B\left(\frac{\Delta \tau}{2}\right) A\left(\frac{\Delta \tau}{2}\right) f^n. \label{eq:Strang}
\end{equation}
So far, we assumed that we have a good approximation of the potential $\phi$ at hand. In practice, this is obtained by solving the QN equation i.e. \eqref{eq:qn} by a spectral Finite Element method (FEM) combined with some Finite Difference (FD) approximations. Since this is not the main part of our work, we refer to the sources above for more details.

In total, the iterative predictor-corrector solution procedure $f^n \rightarrow f^{n+1}$ can be described as follows:
\begin{enumerate}
	\item Given $f^n$, obtain $\phi$ from solving \eqref{eq:qn}.
	\item Given $\phi$, obtain $f^{n+1/2}$ by applying \eqref{eq:Lie}. (predictor step)
	\item Given $f^{n+1/2}$, obtain $\phi$ from solving \eqref{eq:qn} again.
	\item Given $\phi$, obtain $f^{n+1}$ by applying \eqref{eq:Strang}. (corrector step)
\end{enumerate}

This approach has shown to be quite successful, which inter alia is due to the unconditional stability of the SL scheme. Nonetheless, the method lacks structure preserving properties, like conservation of energy and $L^2$-norm. This gives rise to the main idea of this project, which exchanges the non-restrictive time-stepping of the SL method for the structure preservation of an AKW FD scheme, at least for the slow-time subsystem, where a time-step restriction is supposed to not be too computationally expensive. In other words, we solve the poloidal advection equation i.e. \eqref{eq:adv_poloidal} by applying the AKW scheme combined with a suitable time-integration. Since this advection equation essentially consists of a Poisson bracket, which is rich in geometric structure - exactly what the AKW scheme was designed for - we aim at improving the overall conservation properties of the full simulation.

\subsection{Semi-Lagrangian Scheme for the Fast Time Subsystem}
Recapitulating the semi-Lagrangian method, mostly referring to the overview given in \cite{campospinto}, we will start from the general formulation, where the goal is to solve an advection problem of the form
\begin{equation}
    \partial_t f(t,\mathbf{x})+v(t,\mathbf{x})\cdot \nabla f(t,\mathbf{x})=0, \qquad t\in[0,T], \quad \mathbf{x}\in\mathbb{R}^d, \label{eq:adv_toy}
\end{equation}
where $v$ is a velocity field $\mathbb{R}^d\longrightarrow\mathbb{R}^d$, $T > 0$ is the final time and initial conditions are given by $f_0(\mathbf{x})=f(t = 0,\mathbf{x})$. Assuming that $v$ is a given and smooth vector-field, we can use the method of characteristics, i.e. obtaining trajectories $X(t)=X(t;s,x)$ that are solutions to the ODE
\begin{equation}
    X'(t)=v(t,X(t)), \qquad X(s)=x, \qquad t\in[0,T],
\end{equation}
for $\bx \in \bR^d$ and $s \in [0,T]$. It can be shown that the flow $F_{s,t}:x\longrightarrow X(t)$ is invertible and satisfies $(F_{s,t})^{-1}=F_{t,s}$. Thus, we can derive the analytical solution to equation \eqref{eq:adv_toy} as
\begin{equation}
    f(t,\mathbf{x})=f_0((F_{0,t})^{-1}(\bx)),\qquad  t\in[0,T], \quad \mathbf{x}\in\mathbb{R}^d.
\end{equation}
This implies, given two consecutive time-steps $t_n$ and $t_{n+1}$, we can define the backwards flow
\begin{equation}
    B^{n,n+1}=(F_{t_n,t_{n+1}})^{-1},
\end{equation}
in order to advance the solution $f^n$ from time $t_n$ to time $t_{n+1}$, i.e. $f^{n+1} = f^n \circ B^{n, n+1}$. So far, the derivation has been completely analytical. In practice however, we have to introduce approximation errors for discretizing the distribution $f$, e.g. a Spline interpolation on a grid, and for the backwards flow $B$, which generally depends on the discretization of the vector-field $v$.\\

Now we are in the position to apply this methodology to the flux surface and $\vp$-advection, which is, again, discussed in \cite{emily} and \cite{Latu_2017}. For the flux surface advection equation \eqref{eq:adv_flux}, i.e.
\begin{equation}
 \partial_t f + v_\parallel \nabla_\parallel f = 0,
\end{equation}
where we remind that $\nabla_\parallel = \bb \cdot \nabla$. This is a one-dimensional constant velocity advection with velocity $\vp \bb$ on the flux surface $(\theta, z)$ for each given $r$. Thus, we can construct an analytical two-dimensional SL operator, that uses the exact trajectory as the velocity is not related to the flux surface. On the other hand, the $v$-parallel advection operator defined by equation \eqref{eq:adv_par}, i.e.
\begin{equation}
    \partial_t f + \nabla_\parallel \phi\,\, \partial_{v_{\parallel}} f = 0,
\end{equation}
 contains the parallel gradient of $\phi$ which only depends on the spatial coordinates and is therefore constant along $v_\parallel$. As a result, the trajectory used by this one-dimensional SL method can be accurately defined, while the parallel gradient of $\phi$ is computed using a field-aligned FD method. \\
 
When implementing the method, we work on a four dimensional computational grid on which the point values of the distribution functions and potentials for the different time-stages are known or calculated.
Both advection steps use special interpolation techniques, since we will not end up exactly on grid points when tracing back the characteristics, such that they are at least of order three. Additionally, we have to take account for boundary conditions, when the characteristics move outside the computational domain, extending it by extrapolation for instance. \\

\subsection{Arakawa Scheme for the Slow Time Subsystem}
\label{sec:AKW}
Introducing the Arakawa scheme, we mainly reference the original article \cite{Arakawa_1966}, where one is interested in the spatial two-dimensional discretization of the differential equation
\begin{equation}
	\partial_t f + \{\phi, f\} = 0, \label{eq:br_ode}
\end{equation}
where $\phi$ is a given potential and $\{\phi,f\}$ is a Poisson bracket of the form
\begin{equation}\label{eq:poisson_bracket}
	\{\phi,f\} = -  \left(\partial_x\phi\right) \left(\partial_y f\right) + \left(\partial_y\phi\right) \left(\partial_x f\right) \, .
\end{equation}

The main aim of this scheme is the conservation of the following properties
\begin{subequations}\label{conservation-properties}
	\begin{align}
		\text{Mass} : && \dt \int f(t) \d x \text{d} y & = 0 & \Leftrightarrow && \int\bracket{\phi}{f} \d x \text{d} y & = 0 \, , \label{eq:consv_1}\\
		L^2\text{-norm :} && \dt \int f^2(t) \d x \text{d} y & = 0 & \Leftrightarrow && \int f \, \bracket{\phi}{f} \d x \text{d} y & = 0 \, , \\
		\text{Total energy :} && \dt \int \phi \, f(t) \d x \text{d} y & = 0 & \Leftrightarrow && \int \phi \bracket{\phi}{f} \d x \text{d} y & = 0 \, ,
	\end{align}
\end{subequations}
which is deeply embedded in its construction. Albeit we are interested in the application of the scheme in polar coordinates, i.e. a change of coordinates in the Poisson bracket c.f. equation \eqref{eq:bracket_in_polar}, we will for simplicity start with the construction in Cartesian coordinates since conceptionally it does make no difference.

\subsubsection{Construction of the Discrete Bracket}
\label{sec:const_stenc}
\begin{figure}
	\centering
	\begin{minipage}[h]{0.45\textwidth}
		\centering
		\begin{tikzpicture}
			\fill[yellow!40!white] (-1.5,-1.5) rectangle (1.5,1.5);
			
			\draw[step=1.5cm, black, very thin] (-2.4,-2.4) grid (2.4,2.4);
			
			\draw[step=1.5cm, black,very thick] (-1.5001,-1.5001) grid (1.5001,1.5001);
			
			\filldraw [orange]  (1.5,0) circle (2pt);
			\filldraw [orange] (0,1.5) circle (2pt);
			\filldraw [blue] (1.5,1.5) circle (2pt);
			\filldraw [blue] (1.5,-1.5) circle (2pt);
			\filldraw [blue] (-1.5,1.5) circle (2pt);
			\filldraw [orange] (-1.5,0) circle (2pt);
			\filldraw [orange] (0,-1.5) circle (2pt);
			\filldraw [blue] (-1.5,-1.5) circle (2pt);
			\draw[orange] (0, -1.4) -- (0,1.4);
			\draw[orange] (-1.4, 0) -- (1.4, 0);
			\draw[blue, dashed] (1.4, -1.4) -- (-1.4,1.4);
			\draw[blue, dashed] (1.4, 1.4) -- (-1.4, -1.4);
		\end{tikzpicture}
	\end{minipage} \begin{minipage}[h]{0.45\textwidth}
		\centering
		\begin{tikzpicture}
			\fill[yellow!40!white] (-1.5,-1.5) rectangle (1.5,1.5);
			
			\draw[step=1.5cm, black, very thin] (-2.9,-2.9) grid (2.9,2.9);
			
			\draw[step=1.5cm, black,very thick] (-1.5001,-1.5001) grid (1.5001,1.5001);
			
			\filldraw [orange] (3,0) circle (2pt) ;
			\filldraw [orange] (0,3) circle (2pt) ;
			\filldraw [blue] (1.5,1.5) circle (2pt) ;
			\filldraw [blue] (1.5,-1.5) circle (2pt) ;
			\filldraw [blue](-1.5,1.5) circle (2pt) ;
			\filldraw [orange] (-3,0) circle (2pt) ;
			\filldraw [orange] (0,-3) circle (2pt) ;
			\filldraw [blue](-1.5,-1.5) circle (2pt) ;
			\draw[orange] (0, -3) -- (0,3);
			\draw[orange] (-3, 0) -- (3, 0);
			\draw[blue, dashed] (1.4, -1.4) -- (-1.4,1.4);
			\draw[blue, dashed] (1.4, 1.4) -- (-1.4, -1.4);
		\end{tikzpicture}
	\end{minipage}
	\caption{The nine (left) and thirteen (right) point stencil, where the $+$-parts are highlighted in solid orange, while the $\times$-parts are in dashed blue. Here, the middle of the cell corresponds to the point on which a stencil is defined.}
	\label{fig:stencils}
\end{figure}
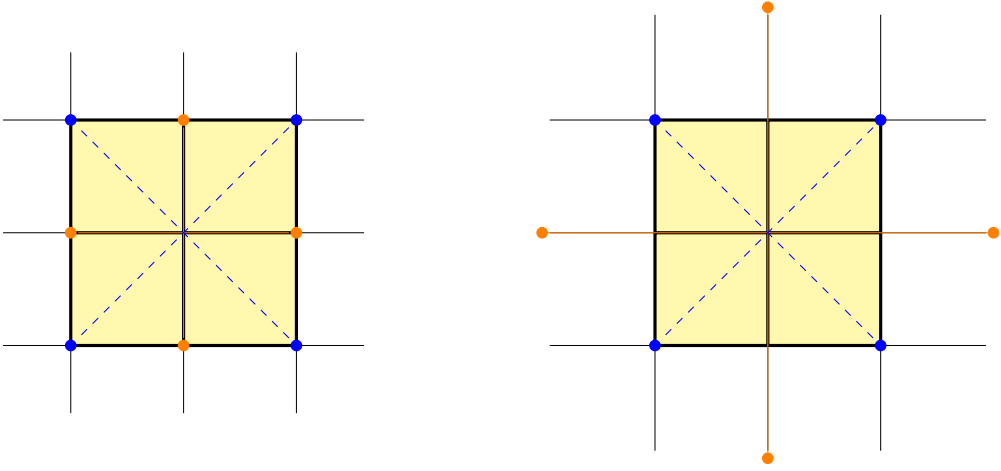

Given a two-dimensional grid $(x_i, y_j)$ for $0 \le i \le N_x, 0 \le j \le N_y$ with a uniform grid size $d >0$, we simplify notation by writing $g_{i,j} = g(x_i, y_j)$ for any function $g$ evaluated at the point $(x_i, y_j)$ and write the collection of point-values as discrete function $g_h$. Denoting the Poisson bracket by 
\begin{equation} \label{eq:J-bracket}
	J(f, g) = \{f, g\}
\end{equation} and following \cite{Arakawa_1966}, we can approximate it at any point $(x_i, y_j)$ as a certain linear combination of the following nine-point stencils, that amounts to
\begin{equation}
	J(f,g)=J_1(f_h, g_h)+\mathcal{O}(d^2), \quad \text{ where } \quad J_1 = \frac{1}{3}(J_1^{++}+J_1^{+\times}+J_1^{\times+}),
\end{equation}
with the stencils defined as
\begin{equation}
	\begin{aligned}
    J_1^{++}& =\frac{1}{4d^2}\left[(f_{i+1,j}-f_{i-1,j})(g_{i,j+1}-g_{i,j-1}) \right. \\
	& \hspace{11mm}  \left.-(f_{i,j+1}-f_{i,j-1})(g_{i+1,j}-g_{i-1,j})\right], \label{eq:stencil1}
	\end{aligned}
\end{equation}
as well as
\begin{equation}
    \begin{aligned}
    	J_1^{+\times} & = \frac{1}{4d^2} \left[ f_{i+1,j} (g_{i+1,j+1} - g_{i+1,j-1}) - f_{i-1,j} (g_{i-1,j+1} - g_{i-1,j-1}) \right. \\
    	& \hspace{11mm} \left. - f_{i,j+1} (g_{i+1,j+1} - g_{i-1,j+1}) + f_{i,j-1}(g_{i+1,j-1} - g_{i-1,j-1}) \right],
    \end{aligned}
\end{equation}
and
\begin{equation}
    \begin{aligned}
    	J_1^{\times+} & = \frac{1}{4d^2} \left[ f_{i+1,j+1} (g_{i,j+1} - g_{i+1,j}) - f_{i-1,j-1} (g_{i-1,j} - g_{i,j-1}) \right. \\
    	& \hspace{11mm} \left. - f_{i-1,j+1} (g_{i,j+1} - g_{i-1,j}) + f_{i+1,j-1} (g_{i+1,j} - g_{i,j-1}) \right].
    \end{aligned}
\end{equation}
This approximation is proven to be the only FD second order approximation of the analytical Poisson bracket that conserves mass, $L^2$-norm and energy for an isotropic mesh. The $+$ and $\times$ notation comes from the patterns done by the FD stencils on the grid, as is visualized in Figure \ref{fig:stencils}. So for example in $J_1^{+\times}$, we choose the $+$-pattern, which corresponds to two central FD schemes, for the discretization of $\partial_x f_{i,j}$ and $\partial_y f_{i,j}$ and the same with the $\times$-pattern for $g_{i,j}$, then multiply the two together getting a discretization of the bracket \eqref{eq:poisson_bracket}.

In order to continue these considerations to fourth order, we introduce $J_2(f_h, g_h)=\frac{1}{3}(J_2^{\times \times} + J_2^{\times+} + J_2^{+\times})$, where we now use the extended twelve point-stencils with the additional
four points $(i+2,j)$, $(i-2,j)$, $(i,j+2)$ and $(i,j-2)$, visualized in Figure \ref{fig:stencils}, such that
\begin{equation}
\begin{aligned}
	J_2^{\times\times} &= \frac{1}{8d^2} \left[(f_{i+1,j+1} - f_{i-1,j-1}) (g_{i-1,j+1} - g_{i+1,j-1}) \right. \\
	& \hspace{11mm} - \left. (f_{i-1,j+1} - f_{i+1,j-1}) (g_{i+1,j+1} - g_{i-1,j-1}) \right] \, ,
\end{aligned}
\end{equation}
as well as
\begin{equation}
	\begin{aligned}
		J_2^{\times+} & = \frac{1}{8d^2} \left[ f_{i+1,j+1} (g_{i,j+2} - g_{i+2,j}) - f_{i-1,j-1} (g_{i-2,j} - g_{i,j-2}) \right. \\
		& \hspace{11mm} \left. - f_{i-1,j+1} (g_{i,j+2} - g_{i-2,j}) + f_{i+1,j-1} (g_{i+2,j} - g_{i,j-2})\right] \, ,
	\end{aligned}
\end{equation}
and
\begin{equation}
	\begin{aligned}
		J_2^{+\times} & = \frac{1}{8d^2} \left[f_{i+2,j} (g_{i+1,j+1} - g_{i+1,j-1}) - f_{i-2,j} (g_{i-1,j+1} - g_{i-1,j-1}) \right. \\
		& \hspace{11mm} \left. - f_{i,j+2} (g_{i+1,j+1} - g_{i-1,j+1}) + f_{i,j-2} (g_{i+1,j-1} - g_{i-1,j-1})\right] \, . \label{eq:stencil2}
	\end{aligned}
\end{equation}
In total, we can now approximate the Poisson bracket $J$ at any point on the grid by
\begin{equation}
	J(f,g) = 2J_1(f_h, g_h)-J_2(f_h, g_h) + \mathcal{O}(d^4),
\end{equation}
up to fourth order while conserving all the desired quantities above as is shown in \cite{Arakawa_1966}. This leads us to define the discrete Poisson bracket $J_h$ on any point of the grid as
\begin{equation} \label{order4-disc-bracket}
	J_{h}(f_h, g_h) = 2J_1(f_h, g_h) - J_2(f_h, g_h),
\end{equation}
thus being able to spatially discretize equation \eqref{eq:br_ode} as
\begin{equation}
	\partial_t f_h = - J_h(\phi_h, f_h),
\end{equation}
such that the right-hand-side is approximated up to order four. By construction, this skew-symmetric discrete bracket satisfies the algebraic properties
\begin{subequations}
	\begin{align}
		\sum_{i,j} J_{h, (i,j)}(\phi_h, f_h) d^2 &= 0,  \\
		\sum_{i,j} f_{i,j} J_{h, (i,j)}(\phi_h, f_h) d^2 &= 0, \\
		\sum_{i,j} \phi_{i,j} J_{h, (i,j)}(\phi_h, f_h) d^2&= 0, 
	\end{align}
\end{subequations}
as is calculated in \cite{Arakawa_1966}, which leads to 
\begin{subequations}
	\begin{align}
		\partial_t \sum_{i,j}  f_{i,j} d^2&= 0,  \\
		\partial_t \sum_{i,j}  f_{i,j}^2 d^2 &= 0, \\
		\partial_t \sum_{i,j}  (\phi f)_{i,j} d^2 &= 0,
	\end{align}
\end{subequations}
where $\partial_t \phi = 0$. This the discrete analogous of the conservation properties of the equations \eqref{conservation-properties}.

\subsubsection{Boundary Conditions} \label{sec:BC}

For points close to or on the boundary of the computational grid, the stencil might use points outside the domain that need to be defined. This is a problem that has to be tackled by introducing boundary conditions (BC) or locally reformulating the stencils. The stencils involve two functions in two spatial variables each, where each direction can be treated individually as we will see shortly and where the motivation comes from the physical properties described in Section \ref{sec:gk-model}. Importantly, the choice of boundary conditions also influences the conservation of mass, $L^2$-norm and energy.

The easiest option is taking periodic BC, as defined in \cite{Arakawa_1966}, where $x_{N_x + 1} = x_1, x_{N_x + 2} = x_2, x_{0} = x_{N_x}, x_{-1} = x_{N_x-1}$, etc. and the same for the $y$-direction. Looking at the stencils \eqref{eq:stencil1}-\eqref{eq:stencil2}, this leaves us only with interior points, where the function-values are known. For our model, this BC is used when looking at point sequences along the $\theta$ direction for both functions as physically this variable is periodic as described in Section \ref{sec:gk-model}.

Secondly, we can define homogeneous Dirichlet BC, i.e. the function values on or outside  the boundary of the grid are equal zero. This means for the stencils \eqref{eq:stencil1}-\eqref{eq:stencil2}, whenever an index is smaller two or greater $N_{x/y}-1$, the corresponding function value is set equal zero. Thus, we end up with reduced stencils, if the point is close to a grid boundary. We will apply this technique to the $r$-direction of the potential $\phi$, which is supposed to be zero outside the interior as described in Section \ref{sec:gk-model}.

Finally, we want to introduce extrapolatory BC; knowing the function values on the boundary and outside the grid, we can directly impose them in the stencils. To make that more precise, we assume $f$ outside the interior of the domain in one variable direction $x$ is known and described by the equilibrium function $f_\text{eq}$. Note that $f_\text{eq}$ still has to fulfil the BC in the other variable direction as we use its values on the boundary. In other words, we directly define the stencils \eqref{eq:stencil1}-\eqref{eq:stencil2}, with the outside values
\begin{equation}
	f_{i, j} = f_\text{eq}(x_i, y_i) \qquad \text{for} \quad  i< 1, i > N_x,  1 \le j \le N_y,
\end{equation}
 where the points in $x$-direction are linearly extended by $d$, e.g. $x_{-1} = x_1 - 2d$. This is exactly the situation of the radial direction of the distribution function $f$ in our model, with an equilibrium function $f_\text{eq}$ that is periodic in $\theta$, c.f. Section \ref{sec:gk-model}.

In \cite{crouseilles2018exponential}, these different kinds of BC with respect to conservation their properties are discussed in more detail, albeit only for the second-order scheme. They show that only the purely periodic BC have perfect conservation properties, for the others one still conserves mass, $L^2$-norm and energy well, but not up to machine precision as there is a small error introduced at the boundary. By their numerical experiments, the AKW scheme works best when combining theses different BC, where they conclude the same combination as we proposed above.

\subsubsection{Transformation to Polar Coordinates}\label{sec:consv-props}
\label{sec:polar}
The original AKW scheme \cite{Arakawa_1966} and all our considerations so far were formulated in Cartesian coordinates on an isotropic mesh, i.e. $\Delta x = \Delta y$. When going to an anisotropic mesh, the convergence order of the stencils still holds true, as we calculate in Appendix \ref{app:Order_of_Arakawa_stencil}. However, the gyro-kinetic model and \texttt{PyGyro} code in \cite{pygyro_code} to which we want to apply the Arakawa scheme, are also formulated in polar coordinates $(r,\theta, z, v_\parallel)$. Transforming the previous part to polar coordinates, we first introduce a polar grid $(r_i, \theta_j)$ for $0 \le i \le N_r, 0 \le j \le N_\theta$ with grid-increments $\Delta r >0$ and $\Delta \theta > 0$. The polar bracket $\{\cdot, \cdot\}^p$ of our model is given by \eqref{eq:bracket_in_polar}, dropping the constant $B_0$, which is different to the 
bracket $\{\cdot, \cdot\}^c$ in Cartesian coordinates from \eqref{eq:poisson_bracket} by a factor of $r$. 
This is due to a change of metric coming from the coordinate transformation from Cartesian to polar coordinates, i.e.
\begin{align}
	x  = r \cos(\theta), \qquad y  = r \sin(\theta),
\end{align}
which gets more clear when looking at the quantities under the integral: 
\begin{subequations}
	\begin{align}
		\int\bracket{\phi}{f}^c \d x \d y &= \int\bracket{\phi}{f}^p r \d r \d \theta, \\
		\Longleftrightarrow   \int  \left[ -\left(\partial_x\phi\right) \left(\partial_y f\right) + \left(\partial_y\phi\right) \left(\partial_x f\right) \right]  \d x \d y &=  \int \left[- \frac{1}{r} \left(\partial_\theta\phi\right) \left(\partial_r f\right) + \frac{1}{r} \left(\partial_r\phi\right) \left(\partial_\theta f\right)\right] r \d r \d \theta, \\
		\Longleftrightarrow  \int J^c(f, \phi)  \d x \d y &= \int J^p(f, \phi) r \d r \d \theta,
	\end{align}
\end{subequations}
where $J^c$ is defined in \eqref{eq:J-bracket} and 
\[ J^p(f,g) = \left[- \frac{1}{r} \left(\partial_\theta\phi\right) \left(\partial_r f\right) + \frac{1}{r} \left(\partial_r\phi\right) \left(\partial_\theta f\right)\right] .\] 
Analogous to the previous Section \ref{sec:const_stenc}, we define the discrete bracket at any point of the polar grid as
\begin{equation}
	J^p_{h, (i,j)}(f_h, g_h) = \frac{1}{r_i} J^c_{h, (i,j)}(f_h, g_h),
\end{equation}
from which we obtain the discrete conserved quantities 
\begin{align} \label{eq:pol_cons_quant}
	 \sum_{i,j}  f_{i,j} r_i \Delta r \Delta \theta = 0,  \qquad \sum_{i,j}  f_{i,j}^2 r_i \Delta r \Delta \theta  = 0, \qquad \sum_{i,j}  (\phi f)_{i,j} r_i \Delta r \Delta \theta  = 0,
\end{align}
for mass, $L^2$-norm and energy respectively. As before, their conservation is equivalent to the equations 
\begin{subequations}\label{eq:alg_prop}
	\begin{align} 
		\sum_{i,j} J^p_{h, (i,j)}(\phi_h, f_h) r_i \Delta r \Delta \theta &= 0,  \\
		\sum_{i,j} f_{i,j} J^p_{h, (i,j)}(\phi_h, f_h) r_i \Delta r \Delta \theta &= 0, \\
		\sum_{i,j} \phi_{i,j} J^p_{h, (i,j)}(\phi_h, f_h) r_i \Delta r \Delta \theta &= 0. 
	\end{align}
\end{subequations}


\section{Numerical Experiments}
\label{sec:num_exp}

While the Arakawa method was implemented in \textit{Python} in order to integrate it seamlessly into the \texttt{PyGyro} code, efforts were made to maximize performance in order to do large scale simulations with more than 300 million degrees of freedom in feasible times. This section is devoted to discuss the implementation of the Arakawa scheme, its integration to the \texttt{PyGyro} code, and further numerical experiments and verifications.

\subsection{Implementation of the Discrete Bracket}

Since the potential $\phi$ is constant while performing the poloidal step, which is defined on a polar domain similar to Section \ref{sec:polar}, and in order to create a computationally efficient scheme, we choose to implement the discrete bracket as a \texttt{scipy} sparse matrix $\mathbb{J}_\phi$ of size $(N_r N_\theta)^2$, that is constructed only dependent on the point-values of $\phi$, mapping the point-values of the current distribution function $f$, such that
\begin{equation}
	\mathbb{J}_\phi f_h = J^p_h(\phi_h, f_h) \approx \{\phi, f\}^p.
\end{equation}
After every time step, we update the non-zero entries of this matrix in-place using the new values of $\phi$. These entries are computed in a \textit{Fortran} routine which was generated using \texttt{pyccel} (see \cite{pyccel}) achieving near-native \textit{Fortran}-performance with much less development time. An explicit time integrator then uses matrix-vector multiplication which is computed in \textit{C} thanks to the usage of \texttt{scipy}. An implicit time stepping makes use of the implemented sparse solvers, also provided by \texttt{scipy}.\\

In order to test if the conserved quantities in equation \eqref{eq:pol_cons_quant} hold, we can directly calculate the equivalent algebraic conditions from equation \eqref{eq:alg_prop}, which should only depend on the definition of $J_h$ and are independent of the actual point values in $f_h$ and $\phi_h$ as long as they satisfy the boundary conditions.  

This is interesting with respect to the discussion of the conservation properties depending on the different BC in \cite{crouseilles2018exponential} and Section \ref{sec:BC}, we therefore implemented all BC discussed in Section \ref{sec:BC}, with additional possible periodic and Dirichlet BC in radial direction of the distribution function $f$. 
\begin{table}[h]
	\begin{tabular}{| l | l | l | l | l | l | l|l|l|}
		\hline
		BC	& Order & Mass & $L^2$-norm & Energy & Order & Mass & $L^2$-norm & Energy \\
		\hline
		Periodic& 2 &  1.47e-14 & 2.62e-14 & 1.50e-14 & 4 & 3.93e-14 & 5.57e-14 & 7.44e-14  \\ \hline
		Dirichlet& 2 &  1.35e-13 & 9.09e-13 & 8.67e-13 & 4 & 4.12e-13 & 4.37e-11 & 3.40e-12 \\ \hline
		Extrapolation& 2 &  8.53e-14 & 1.09e-11 & 4.57e-12 & 4 & 2.56e-13 & 2.55e-11 & 1.63e-11 \\
		\hline
	\end{tabular}
	\medskip
	\caption{Algebraic conservation properties, i.e. equation \eqref{eq:alg_prop}, for vectors of size $N_rN_\theta$, with $N_r = N_\theta = 64$, where $f_h \in \bR^{N_rN_\theta}$ and $\phi_h \in \bR^{N_rN_\theta}$ have uniformly distributed values between $-100$ and $100$, while satisfying the BC.}
	\label{tab:alg_prop}
\end{table}

The results of this first test can be found in Table \ref{tab:alg_prop}. Even though the values are not on machine-precision, as was predicted in \cite{crouseilles2018exponential} and which is due the imperfect conservation at the non-periodic boundaries, they are still very small, and we expect good conservation properties from using this scheme in the full code. We also observed, that these values increase proportionally to the norm of $\phi$, which may should be normalized in the indicator.

\newpage

\subsection{Poloidal Advection Test-Case}

In order to further verify the AKW scheme and its properties, we look at a test-case similar to Section 3.5.2 in \cite{emily}, that is, we consider a poloidal advection on a polar domain with $r \in [1, 20]$ and $\theta \in [0, 2\pi]$. The equilibrium distribution function $f_\text{eq}$ reads
\begin{equation}
	f_\text{eq}(r, v_\parallel) = \frac{n_0(r)}{\sqrt{2\pi T_i(r)}} \exp\left( - \frac{v_\parallel^2}{2 T_i(r)} \right)
\end{equation}
with radial profiles
\begin{equation}
	\mathcal{P}(r) = C_{\mathcal{P}} \exp\left[ - \kappa_\mathcal{P} \delta r_\mathcal{P} \tanh \left( \frac{r - r_p}{\delta r_\mathcal{P}} \right) \right]
\end{equation}
for $\mathcal{P} \in \{T_i, n_0\}$ and with constants
\begin{align}
	C_{T_i} & = 1 & C_{n_0} & = \left(r_\text{max} - r_\text{min}\right) \left[\int_{r_\text{min}}^{r_\text{max}} \exp\left[ - \kappa_{n_0} \delta r_{n_0} \tanh \left( \frac{r - r_p}{\delta r_\mathcal{P}} \right) \right] \d r \right]^{-1}
\end{align}
and parameters
\begin{align*}
	B_0 & = 0 \, , & R_0 & = 239.8081535 \, , & r_\text{min} & = 0.1 \, , & r_\text{max} & = 14.5 \, , & r_p & = \frac{r_\text{max} - r_\text{min}}{2} \, , \\
	\epsilon & = 10^{-6} \, , & \kappa_{n_0} & = 0.055 \, , & \kappa_{T_i} & = 0.27586 \, , & \delta r_{T_i} & = \frac{\delta r_{n_0}}{2} = 1.45 \, , & \delta r & = \frac{4 \delta r_{n_0}}{\delta r_{T_i}} \, .
\end{align*}

Initial potential $\phi$ and initial distribution function $f$ are given by
\begin{subequations}
	\begin{align}
		\phi(r, \theta) & = -5 r^2 + \sin(\theta), \\
		f(t = 0, r, \theta) & = f_\text{eq}(r, \vp=0) + B(r, \theta), \\
		B(r, \theta) & = \begin{cases}
			\cos(\frac{\pi}{8} \sqrt{(r - 7)^2 + 2(\theta - \pi)^2}), & \text{if } \sqrt{(r - 7)^2 + 2(\theta - \pi)^2} \le 4, \\
			0, & \text{else}.
		\end{cases},
	\end{align}
\end{subequations}
where the radius-dependent equilibrium function $f_\text{eq}(r)$ is the same as in \cite{Latu_2017} with a fixed $\vp = 0$. For this system, the exact trajectories are given by
\begin{subequations}
	\begin{align}
		\theta(t) &= \theta_0 - 10 t, \\
		r(t) &= \sqrt{ r_0^2 + \frac15 \left( \sin\left(\theta_0 - 10 t\right) - \sin\left(\theta_0 \right)\right)},
	\end{align}
\end{subequations}
for initial points $(r_0, \theta_0)$. Following these characteristics with the initial distribution yields an exact solution to our problem, that can be used to calculate the error of our numerical scheme. The initial distribution is of Gaussian-like shape, as displayed in Figure \ref{fig:init_f}, and gets advected around the center while keeping its shape.

\begin{center}
	\begin{minipage}[h]{0.5\textwidth}
		\centering
		\begin{figure}[H]
			\includegraphics[width=\linewidth]{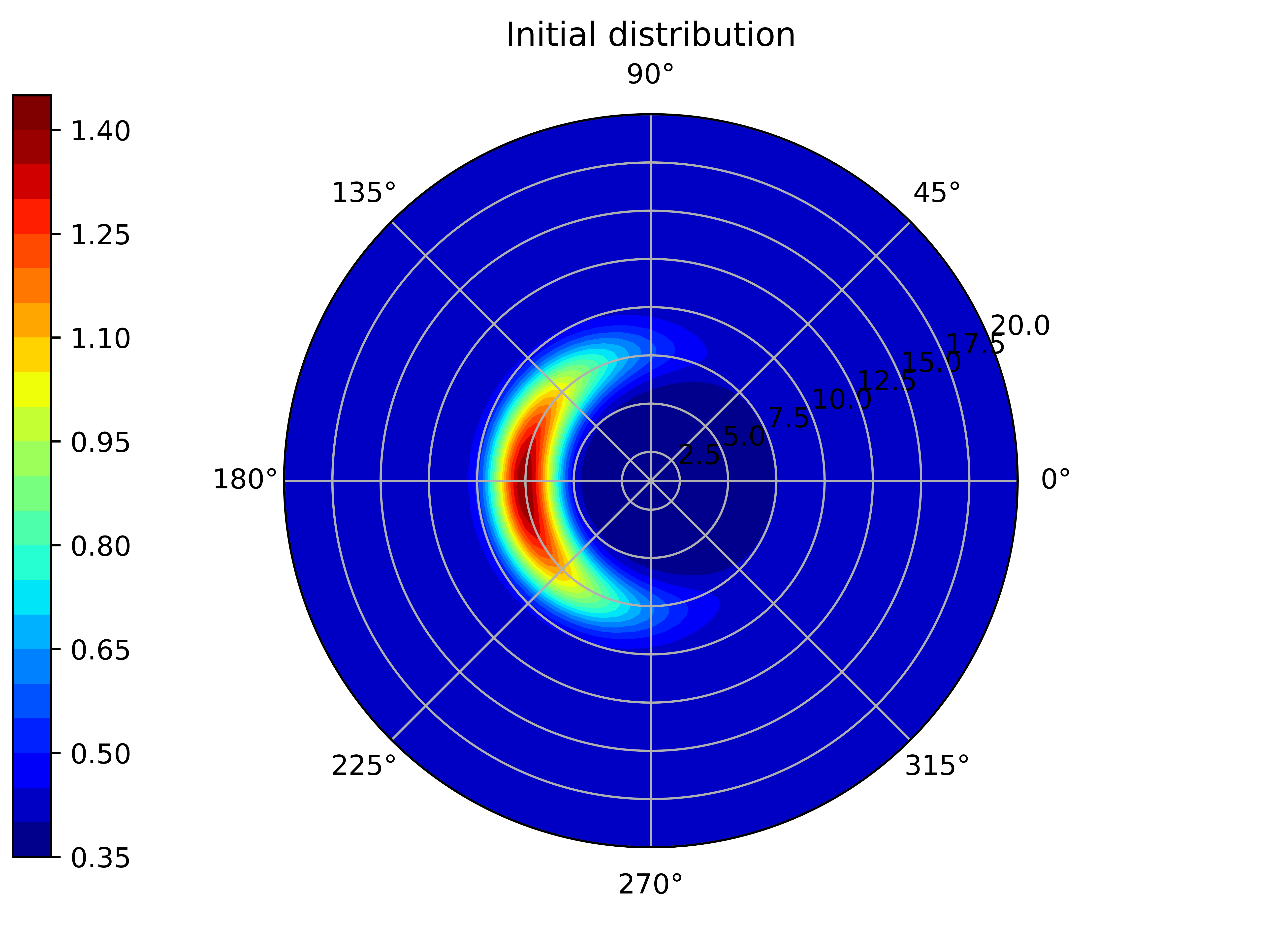}
			\vspace{-0.6cm}
			\caption{The density $f$ at initial time.}
			\label{fig:init_f}
		\end{figure}
	\end{minipage}\begin{minipage}[h]{0.5\textwidth}
		\centering
		\begin{figure}[H]
		\includegraphics[width=\linewidth]{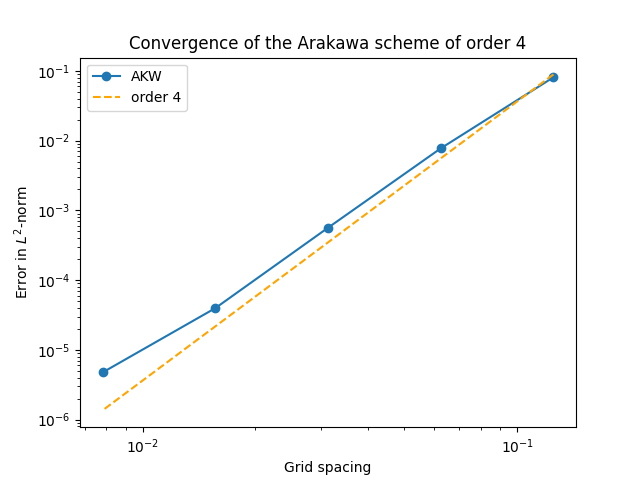}
		\vspace{-0.6cm}
		\caption{Convergence curve of the $L^2$-error from Table \ref{tab:num_exp}.}
		\label{fig:vortex}
	\end{figure}
	\end{minipage}
\end{center}
\vspace{0.3cm}

We use the AKW scheme as introduced in Section \ref{sec:AKW} with extrapolatory BC, given by the equilibrium function, and different mesh sizes. The time integration is done by an explicit Runge-Kutta scheme of order $4$, where we take the step-size $\Delta t = 0.001 / N_r$ and perform a total number of $N = T_\text{end} / dt$ steps to the final time $T_\text{end} = 0.02$. This way, the time-discretization is fine enough, such that the dominant error is coming from the spatial discretization of the bracket only. The numerical results are listed in Table \ref{tab:num_exp}.
\begin{table}[h]
    \begin{tabular}{| l | l | l | l | l | l | l|l|l|}
    \hline
     $N_r$& $N_\theta$ & $L^2$-error & order & conserved mass & conserved $L^2$-norm & conserved energy \\
	\hline
	16 & 16 & 7.72e-03 &  &9.20e-10 & 1.43e-09 & 9.56e-12\\ \hline
	32 & 32 & 5.62e-04 & 3.78 &8.06e-10 & 1.25e-09 & 5.79e-12\\ \hline
	64 & 64 & 3.96e-05 & 3.83 &1.01e-09 & 1.57e-09 & 6.49e-12\\ \hline
	128 & 128 & 4.74e-06 & 3.06 &1.04e-09 & 1.61e-09 & 5.82e-12\\ \hline
    \end{tabular}
	\medskip
	\caption{Spatial discretization and conservation errors for different mesh sizes $(N_r, N_\theta)$. The error of the solution with respect to the exact solution is calculated in the $L^2$-norm. The conservation errors are relative errors of the initial discrete quantities compared to their values at final time.}
	\label{tab:num_exp}
\end{table}
The error is not up to the desired order of $4$, this could be due to conflicting BC as discussed in \cite{crouseilles2018exponential}, but it is not far off and converges towards it for smaller time-steps. The conserved properties look similar to the ones in Table \ref{tab:alg_prop}, albeit being a few orders of magnitude higher which seems to be due to the bigger norm of $\phi$ and imperfect BC.  

\begin{figure}[h]
	\includegraphics[width=0.45\linewidth]{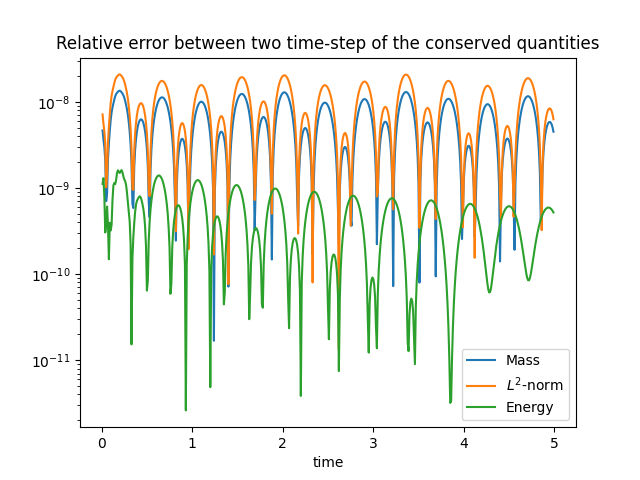}
	\includegraphics[width=0.45\linewidth]{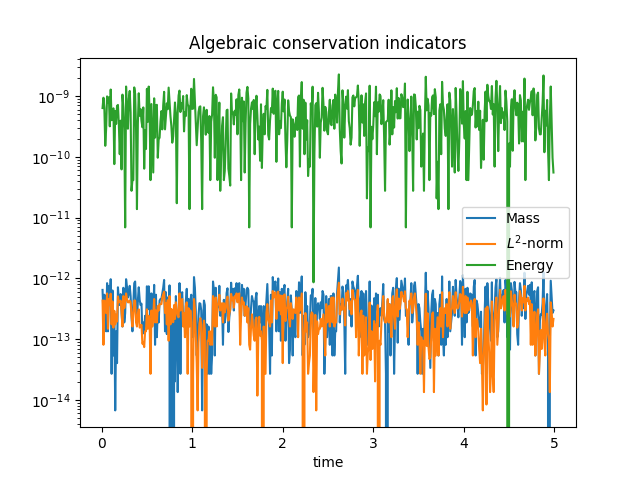}
	\vspace{-0.3cm}
	\caption{The relative error between two time-steps of the conserved quantities (left), i.e. equation \eqref{eq:pol_cons_quant},  and their algebraic indicators at that time-step (right), i.e. equation \eqref{eq:alg_prop}.}
	\label{fig:cons}
\end{figure}

As another test, we solve the system over a longer period of time and keep track of the mass, the $L^2$-norm and the energy. Taking $N_r = N_\theta = 64$, time-step $\Delta t = 0.01$ and perform $N = 500$ steps, yields relative errors of the conserved quantities and values of algebraic indicators as shown in Figure \ref{fig:cons}. 

The long time conservation seems to be in accordance with the results in Tables \ref{tab:num_exp} and \ref{tab:alg_prop}, albeit the algebraic indicators being larger than before due to the larger norm of $\phi$, as was observed before.

\subsection{Full Gyro-kinetic Simulations}

The \texttt{PyGyro} code is a \textit{Python 3} library for gyro-kinetic simulations leveraging the acceleration provided by the modules \texttt{Pythran} (see \cite{Pythran}), \texttt{Numba} (see \cite{Numba}), or \texttt{Pyccel} (see \cite{pyccel}). It is highly parallelized using \texttt{MPI} and thus suitable for running even large-scale simulations on computing clusters. In the following, we will look at results from simulations of the full gyro-kinetic model of \cite{emily}, and then compare the new combination of Semi-Lagrangian and the Arakawa method to the results purely using the Semi-Lagrangian scheme. We will pay special attention to the energy conservation in the poloidal step. Since turbulences mostly occur in this substep, this is also where preserving energy is most crucial.

The parameters used are the same as those in the previous section. The grid size is as follows:
\begin{align*}
	& N_r = 128, && N_\theta = 256, &&& N_z = 128, &&&& N_{v_\parallel} = 72.
\end{align*}
The simulation also uses the following parameters for the initial perturbation:
\begin{align*}
	\iota & = 0, & m & \in \left\{5, 15\right\}, & n & = 1.
\end{align*}
The equilibrium function $f_\text{eq}$ as well as the radial profiles $\{T_i , T_e, n_0\}$ are defined in in the previous section.

\begin{figure}
	\centering
	\includegraphics[width=\linewidth]{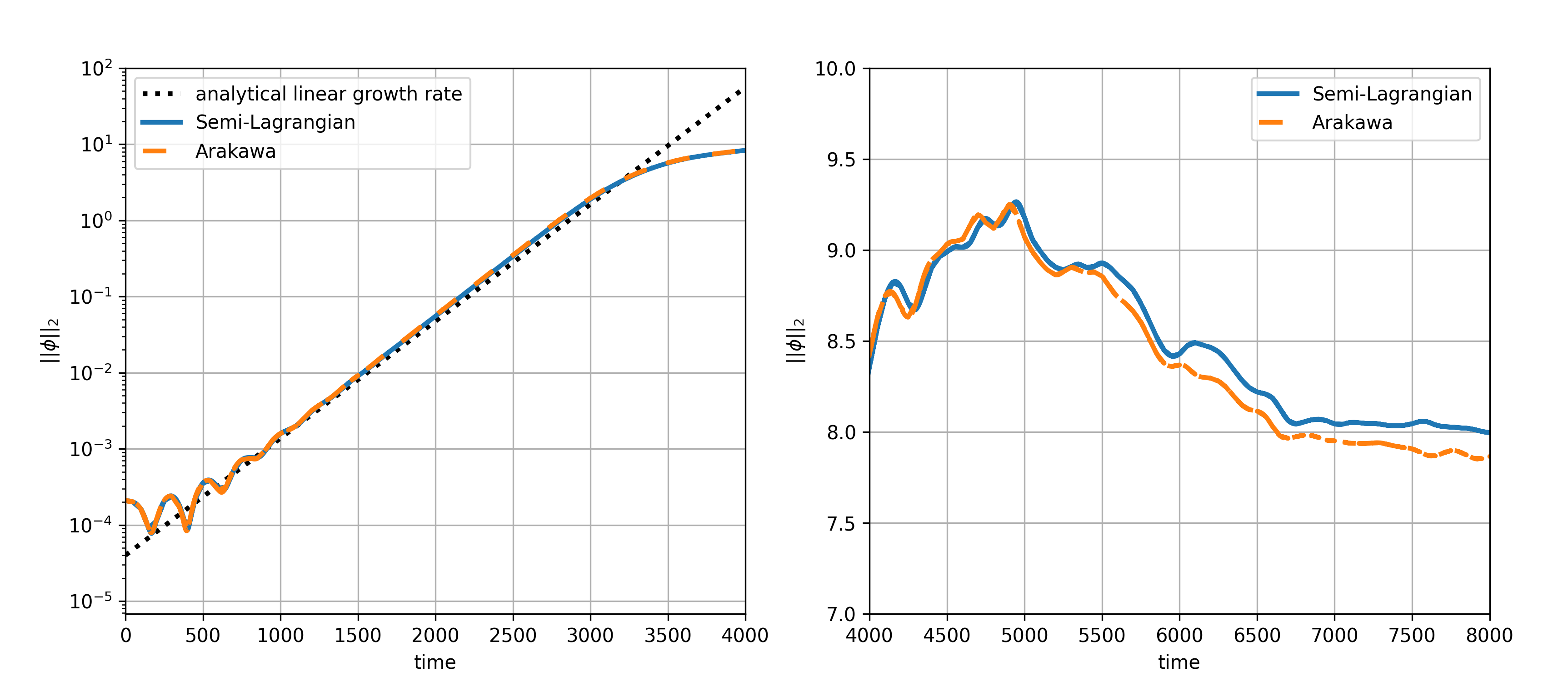}
	\vspace{-0.5cm}
	\caption{The $L^2$-norm of the electric potential $\phi$ comparing the Arakawa method and the Semi-Lagrangian scheme. The results of both simulations closely follow the analytical growth rate of $\norm{\phi}_2 = 4\cdot 10^{-5} \times \exp\left(0.00354 t\right)$ in the linear regime (until $t\sim 3000$) and then saturate in the order of magnitude 10. Observe that the left plot is with a semi-logarithmic scale on the $y$-axis, while the right plot is linear on both axes.}
	\label{fig:l2phi}
\end{figure}

Firstly we compare the Arakawa method to the second-order Semi-Lagrangian scheme by looking at the $L^2$-norm of the electric potential $\phi$, c.f. Figure \ref{fig:l2phi}. We observe the expected behaviour for both schemes where $\norm{\phi}_2$ follows the analytical growth rate of $\norm{\phi}_2 = 4\cdot 10^{-5} \times \exp\left(0.00354 t\right)$ very well in the linear regime (until $t \sim 3000$) (shown on the left in a semi-logarithmic plot). After that the $L^2$-norm saturates and settles at a value of around $8.0$. Notably, the two simulations agree very well even in the non-linear regime which is a good consistency check.

\begin{figure}[h]
	\centering\hspace{-0.5cm}
	\includegraphics[width=0.46\linewidth]{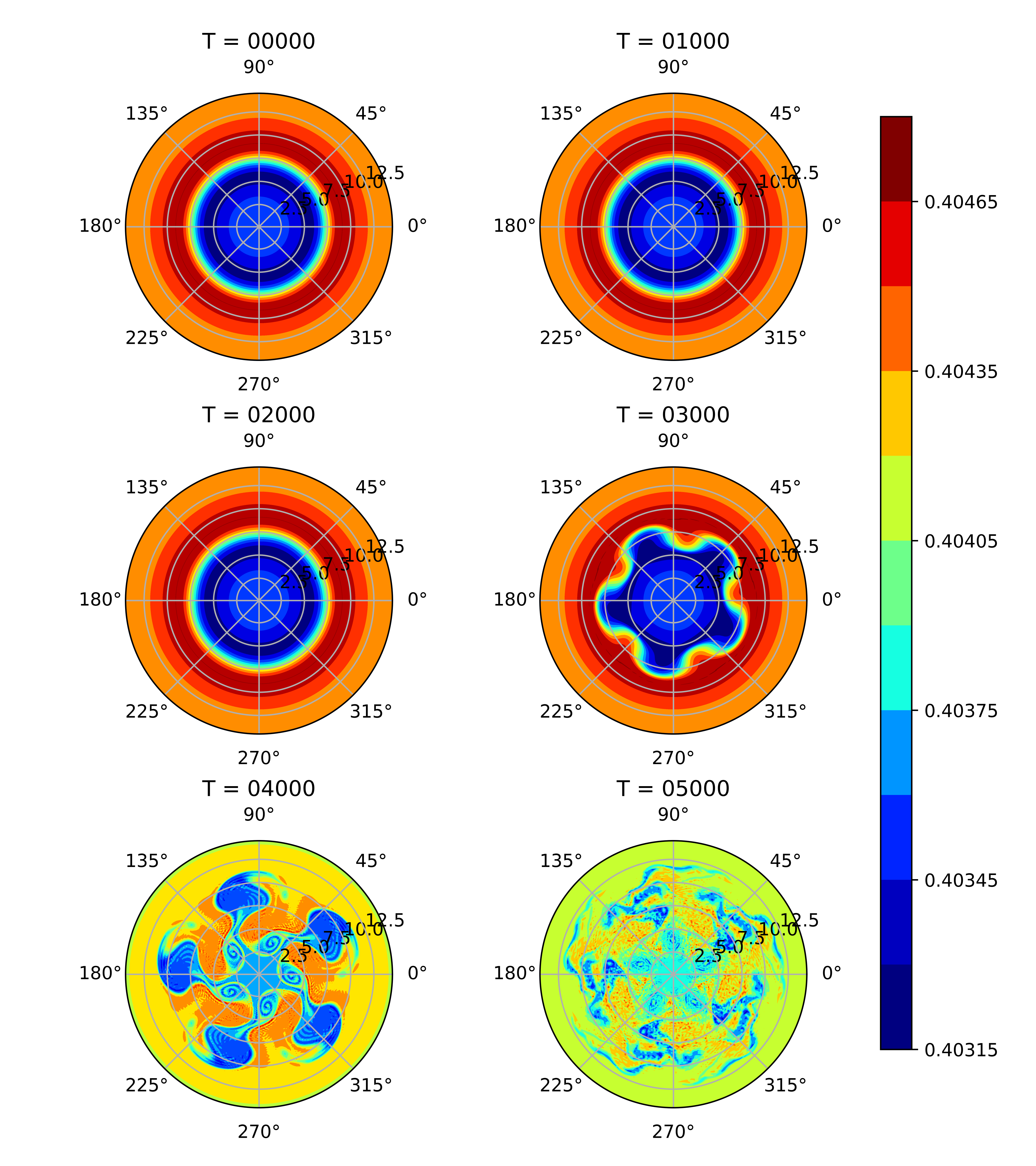}
	\includegraphics[width=0.46\linewidth]{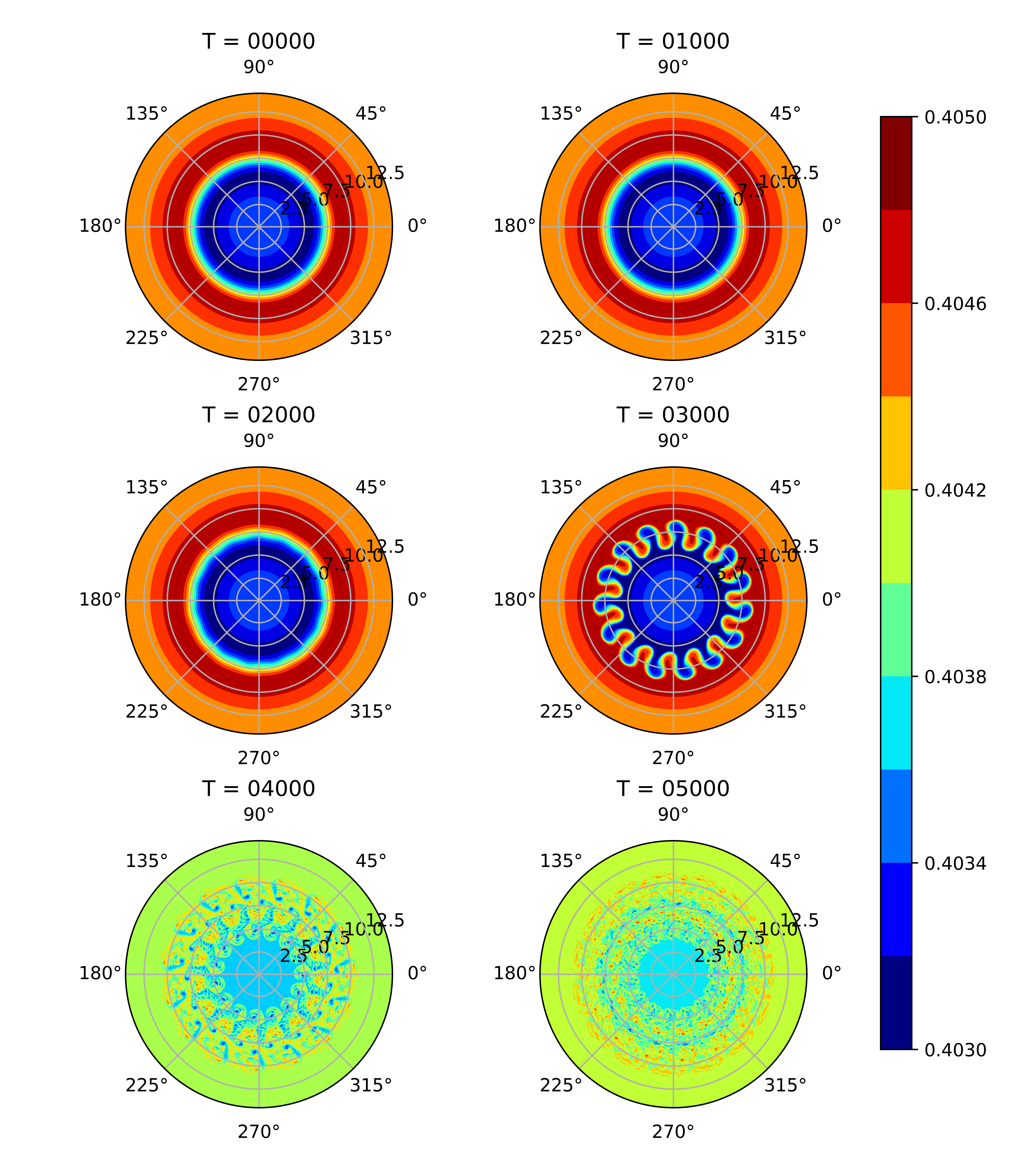} \\
\hspace{-0.5cm}
	\includegraphics[width=0.47\linewidth]{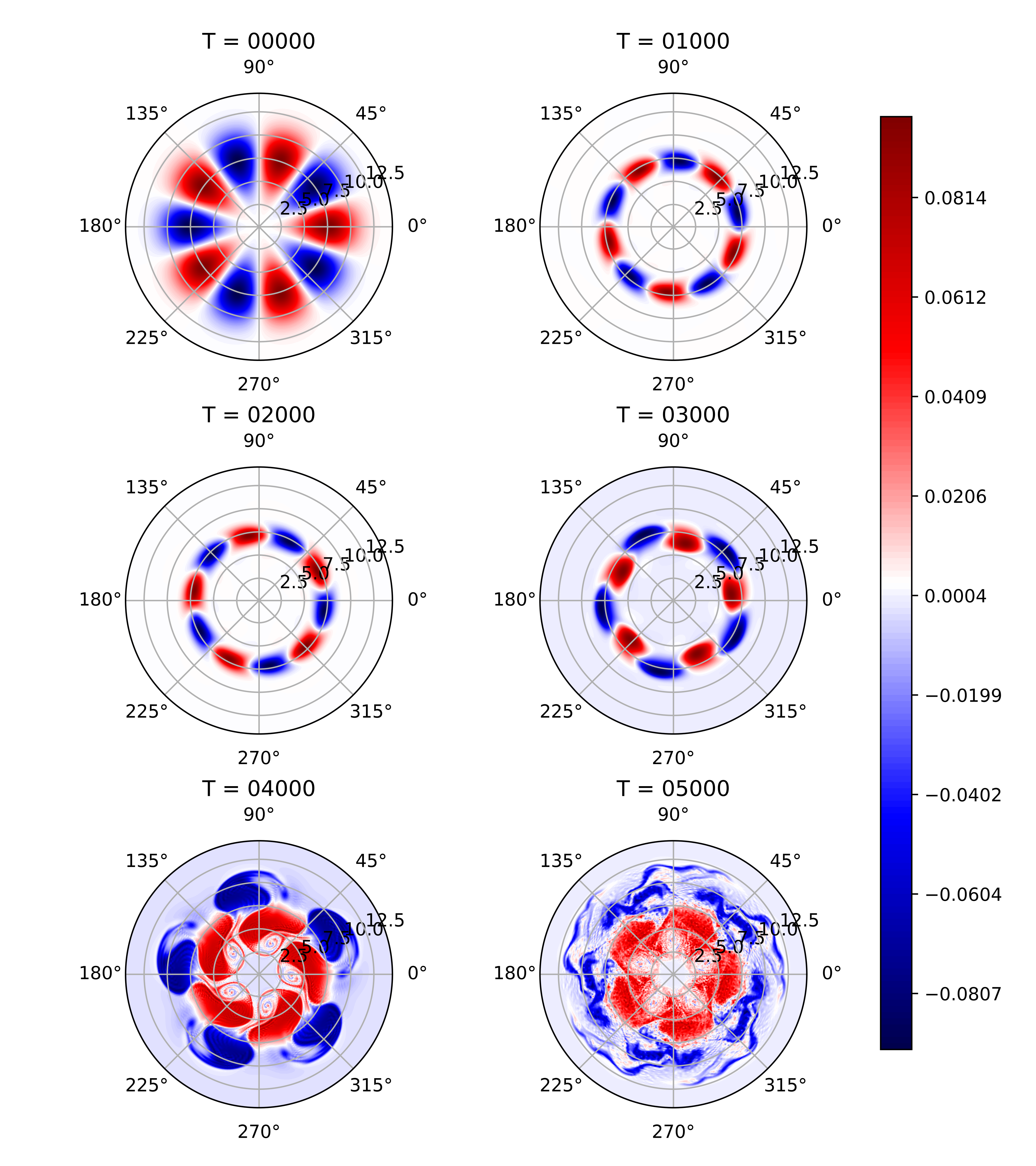}
	\hspace{-0.4cm}
	\includegraphics[width=0.47\linewidth]{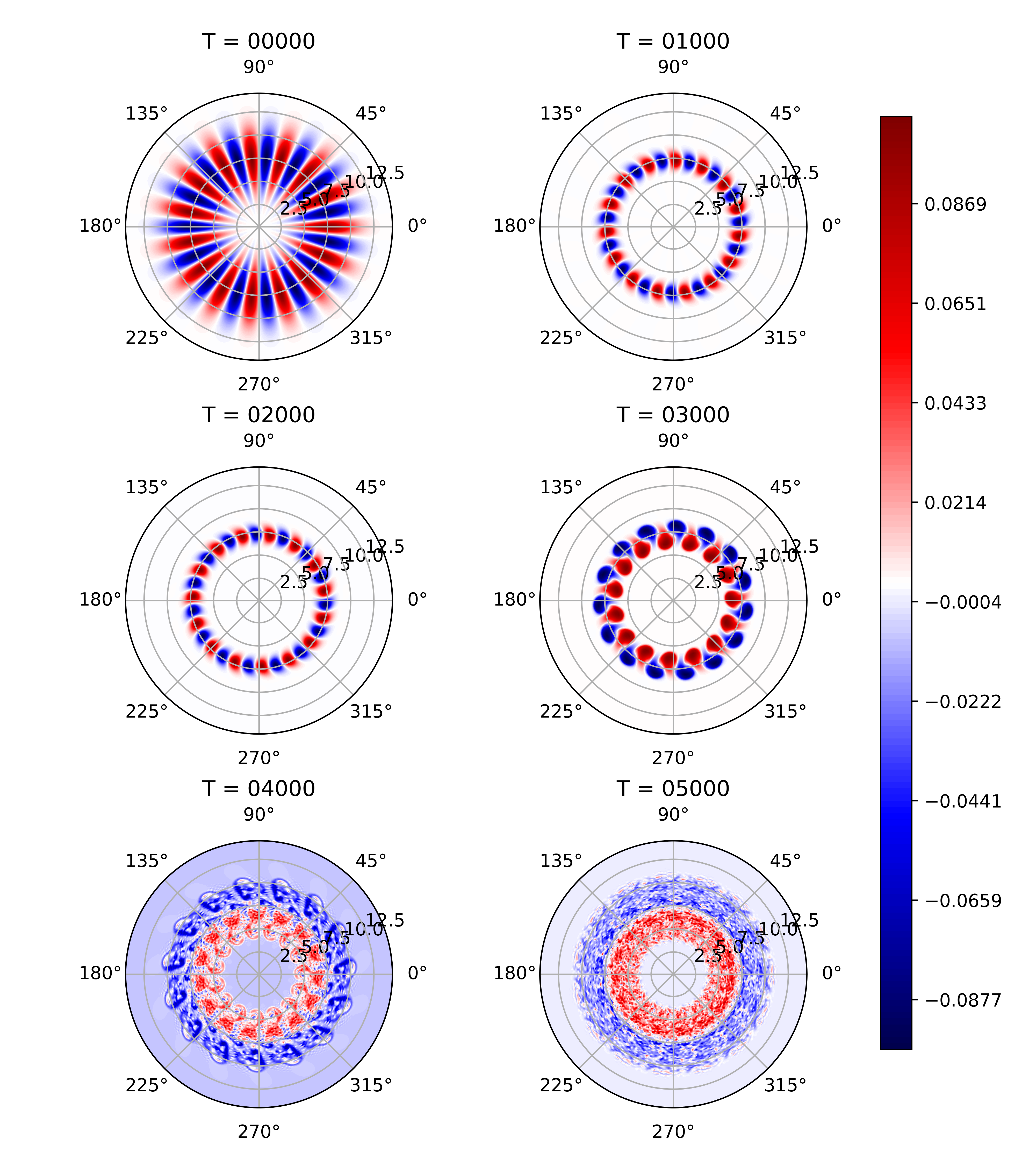}
	\vspace{-0.3cm}
	\caption{The distribution function (above) and the difference between the distribution function and the equilibrium distribution function (below) for different points in time for $m=5$ (left) and $m=15$ (right). Shown are $(r,\theta)$-slices at $(z,v_\parallel) = (0,0)$ for a simulation with the Arakawa scheme.}
	\label{fig:akwfullndiff}
\end{figure}

Furthermore, we plot both the full distribution function and additionally its difference to the equilibrium distribution function in Figure \ref{fig:akwfullndiff}, up to $t=5000$ in steps of $1000$. Shown are slices at $z=0$ and $v_\parallel \simeq 0$ in a polar projection with variables $r$ and $\theta$. We see that it behaves as expected with turbulences forming after the linear regime is over, appearing in the plot for $t=4000$. Interesting to note is the fact that the difference to the equilibrium function seems to increase even for low-valued regions of the domain, although only by a relative difference of half a percentage point.
This effect also occurs in the simulation using the SL scheme, to an even larger extend, so it is not an artifact of the Arakawa scheme. 
One point of discussion is, if this error is due to the boundary conditions, but compared to extrapolation BC, we see no improvement when switching to Dirichlet BC, in this case constant extrapolation by the boundary value, which seem to be the only other plausible option.



Our main point of comparison between the Arakawa method and the Semi-Lagrangian scheme is the conserved quantities in \eqref{conservation-properties}: the mass and $L^2$-norm of the distribution function, and the potential energy, as well as the kinetic energy
\begin{equation}
	E_\text{kin} = \frac{m}{2} \int \left(f(t, r, \theta, z, v_\parallel) - f_\text{eq}(r, v_\parallel)\right) v^2_\parallel \d r \d \theta \d z \d v_\parallel \, .
\end{equation} \\

\begin{figure}[h]
		\centering
		\includegraphics[width=0.8\linewidth]{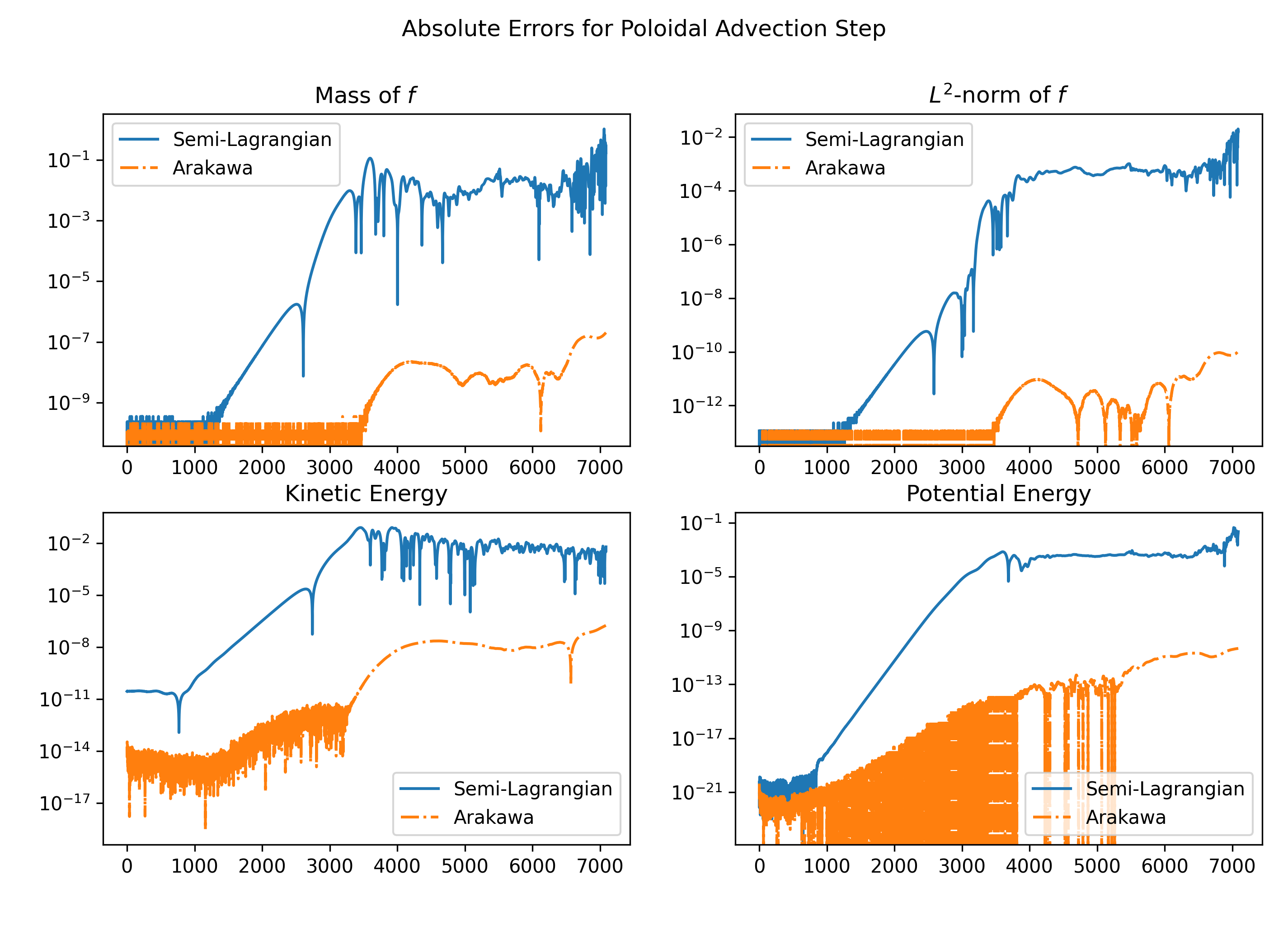}\\
		\vspace{-0.2cm} \hspace{0.3cm}
		\includegraphics[width=0.78\linewidth]{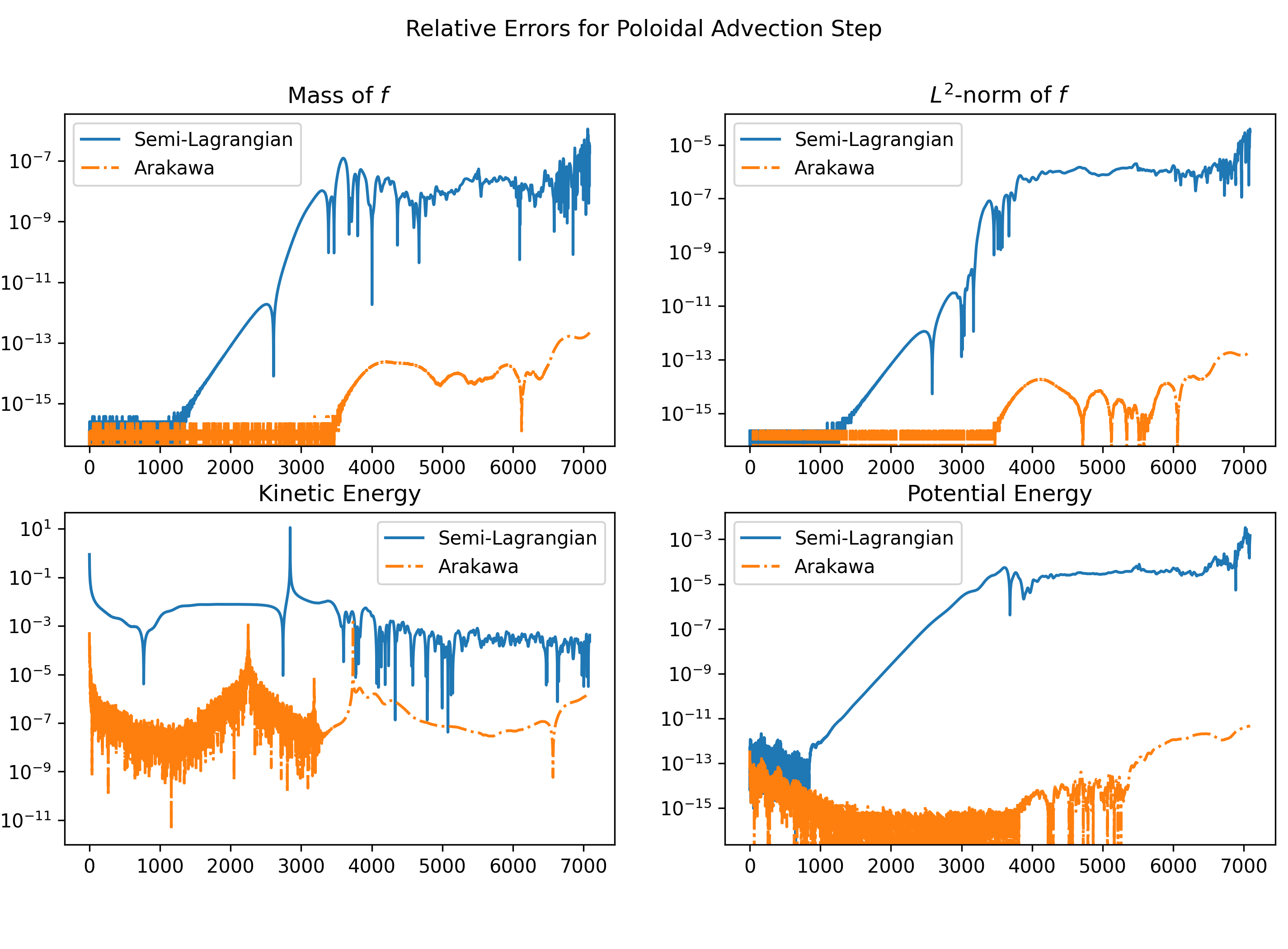}
		\vspace{-0.8cm}
		\caption{Conservation errors corresponding to the poloidal advection step: The quantities are computed once before the poloidal advection step, and then afterwards, and their difference is plotted. The simulation was done with $m=15$. Shown are the absolute (above) and relative (below) errors on a semi-logarithmic scale. It can be seen that the Arakawa method improves significantly the conservation of all four quantities. In particular, they are conserved up to machine precision in the linear regime (until $t\sim 3000$).}
		\label{fig:absnrelerrlog}
		\vspace{-2cm}
	\end{figure}

From a simulation with the above grid size and time step-size $\Delta t = 1$, we investigate the conservation properties for the four above mentioned quantities by computing them before and after the poloidal advection step and plotting the absolute and relative errors in Figure \ref{fig:absnrelerrlog} on a semi-logarithmic scale.
\newpage

\newpage

\newpage

We can clearly see that the Arakawa scheme preserves the conserved quantities much better than the Semi-Lagrangian scheme, although the error is of order of machine precision only in the linear phase. After that, the error grows but remains multiple orders of magnitude smaller than for the semi-Lagrangian scheme. The kinetic energy is also much better preserved which was to be expected, since the behaviour of conservation acts similarly to the mass of the distribution function but with $\vp^2$ weights on each $\vp$-slice, keeping in mind that $\vp$ is just a parameter in each poloidal advection step. The additional error that is introduced by the imperfect boundary conditions is responsible for the conservation properties being bigger than the machine precision; this is especially true for later times when the distribution function moves further away from the equilibrium that is assumed to hold outside the domain. \\

\begin{figure}[h]
	\centering
	\includegraphics[width=0.9\linewidth]{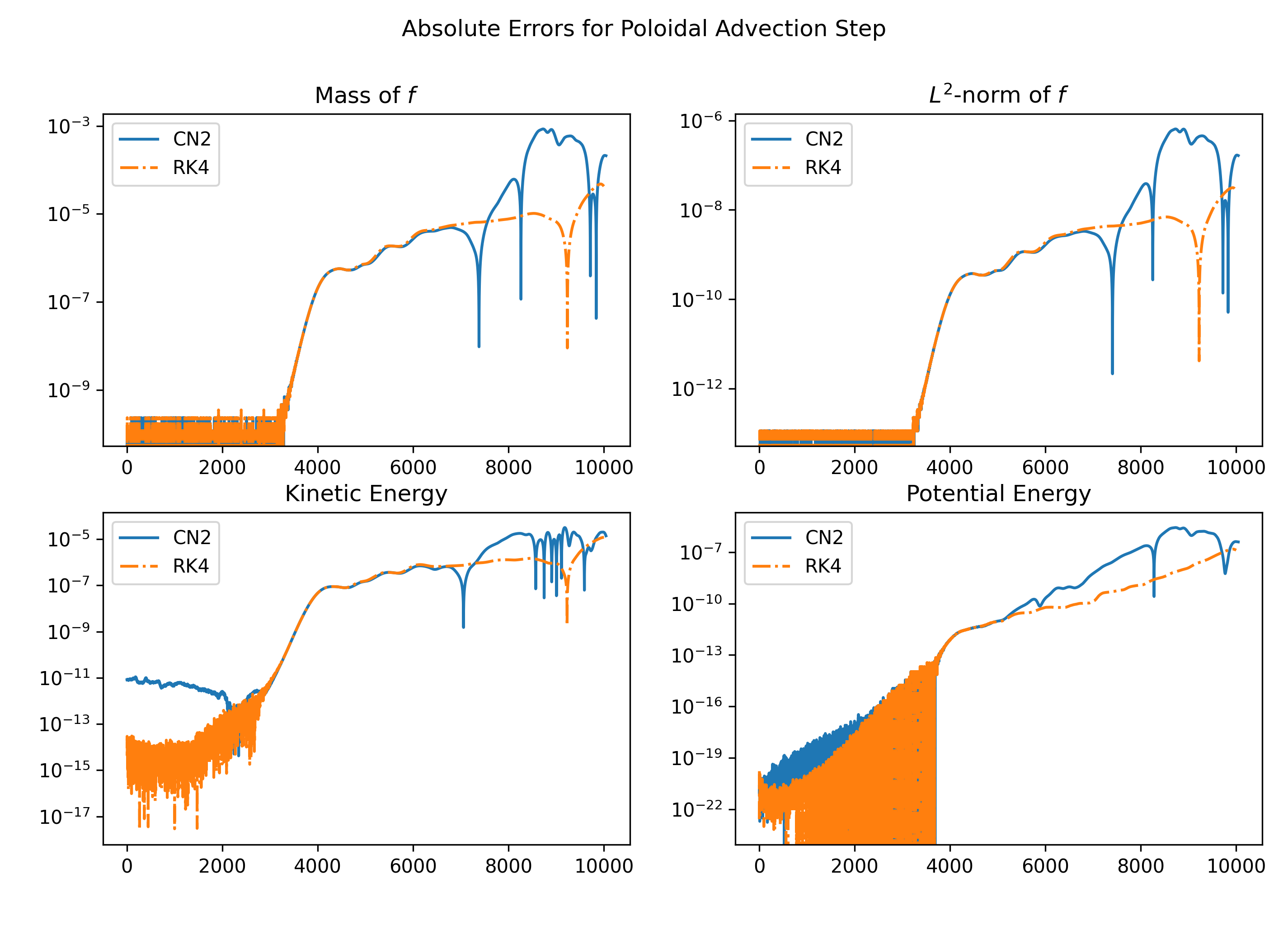}
	\vspace{-0.8cm}
	\caption{Shown are the absolute errors that are computed as described in Figure \ref{fig:absnrelerrlog} on a semi-logarithmic scale, comparing the implicit Crank-Nicolson time integrator, which is of order 2, to the explicit Runge-Kutte scheme of order 4. The errors are quite similar for both methods.}
	\label{fig:abserrlogakw}
\end{figure}

Lastly, we compare the two time-integrators in Figure \ref{fig:abserrlogakw} as described in the text: The second order Crank-Nicolson method and the fourth-order Runge-Kutta scheme, both for the full model simulation. For the RK4 scheme we have a CFL condition which reads
\begin{equation}
	CFL = \max\left(\mathbb{J}_\phi\right) \frac{\Dt}{\Delta x}
\end{equation}
This becomes restrictive only in the non-linear phase later on in the simulation. One can thus use the explicit RK4 with large time step-size $\Dt$, and either use the implicit CN2 method or the RK4 for the later time domain. By observing the accuracy in the discrete conservation of continuous invariants, we see that the difference between the two is negligible for early times, but the RK4 performs better for the very late phase.

\newpage


\section{Conclusion}
\label{sec:conclusion}

In this work we studied a splitting for the gyro-kinetic equation in fast and slow subsystems where demonstrated that the Arakawa scheme yields superior conservation properties when implemented as the spatial discretization of the Poisson bracket defining the poloidal split-step, compared to a full Semi-Lagrangian scheme. This conservation concerns the mass and $L^2$-norm of the distribution function as well as the physical energies $E_\text{kin}$ and $E_\text{pot}$. The improvement was measured by computing the error occurring during the poloidal advection substep of the system. Although the mass and $L^2$-norm of the distribution function, as well as the potential energy are preserved up to machine precision only in the linear phase, the error is for all the quantities multiple orders of magnitude lower compared to the Semi-Lagrangian implementation. The usage of the more physical boundary conditions of extrapolating the distribution function with its equilibrium values outside the domain, has shown to be the right choice, albeit introducing a non-negligible error at later times. Manipulating the Arakawa scheme, to avoid imperfect boundary conditions and fall better in the framework of this problem, is left for future work. 
The comparison between the time-integrators, namely the second order Crank-Nicolson scheme and the 4th order Runge-Kutta method, showed that their behaviour is identical except for very late times, where strong turbulences occur.

\begin{acknowledgement}
	\textbf{Acknowledgements}\\
	The authors would like to thank Emmanuel Franck, Hélène Hivert, Guillaume Latu, Hélène Leman, Bertrand Maury, Michel Mehrenberger, and Laurent Navoret for organizing the CEMRACS conference 2022 and for the wonderful opportunity to come to Marseille and do research. We express special thanks to Michel Mehrenberger and Virginie Grandgirard for the daily supervision and general shaping of the project, and to Xue Hong for discussions on the theoretical side. Dominik Bell and Frederik Schnack thank Eric Sonnendrücker for the opportunity to participate in the CEMRACS and the fruitful discussions after the conference to give this project the final details. They also thank Pierre Navaro for great discussions on and off topic.
\end{acknowledgement}

\bibliographystyle{ieeetr}
\bibliography{literature}

\appendix
\addtocontents{toc}{\protect\setcounter{tocdepth}{0}}

\section{Order of the Arakawa Stencils}\label{app:Order_of_Arakawa_stencil}

For isotropic grids, the approximation properties of the Arakawa scheme has been shown in \cite{Arakawa_1966}. The following \textit{Mathematica} code calculates the approximation power of the discrete stencils on anisotropic meshes. First, we define the series expansion of $f$ and $\phi$ up to order 4:
\begin{subequations}
	\begin{align}
		\mathrm{F}(i, j) & = \text{Series}[f(x + i \Delta x , y+ j \Delta y),\{\Delta x,0,3\},\{\Delta y,0,3 \}] \,, \\
		\mathrm{Phi}(i, j) & = \text{Series}[\phi(x + i \Delta x , y+ j \Delta y),\{\Delta x,0,3\},\{\Delta y,0,3\}] \,.
	\end{align}
\end{subequations}

\subsubsection*{Second Order Scheme (Nine-Point-Stencils)}

Implementing the stencils from section \ref{sec:const_stenc}:
\begin{subequations}
	\begin{align}
		& \begin{aligned}
			\text{J1pp} & = \frac{1}{4 \Delta x \Delta y} \left[ (\text{F}(1,0)-\text{F}(-1,0)) (\text{Phi}(0,1)-\text{Phi}(0,-1)) \right. \\
		& \qquad \hphantom{\frac{1}{4 \Delta x \Delta y}} \left. - (\text{F}(0,1)-\text{F}(0,-1)) (\text{Phi}(1,0)-\text{Phi}(-1,0)) \right] ,
	    \end{aligned} \\
		& \begin{aligned}
			\text{J1px} & = \frac{1}{4 \Delta x \Delta y} \left[ -\text{F}(-1,0) (\text{Phi}(-1,1) - \text{Phi}(-1,-1)) + \text{F}(0,-1) (\text{Phi}(1,-1)-\text{Phi}(-1,-1)) \right. \\
			& \qquad \hphantom{\frac{1}{4 \Delta x \Delta y}} \left. -\text{F}(0,1) (\text{Phi}(1,1)-\text{Phi}(-1,1))+\text{F}(1,0) (\text{Phi}(1,1)-\text{Phi}(1,-1))\right],
		\end{aligned} \\
		& \begin{aligned}
			\text{J1xp} & = \frac{1}{4 \Delta x\Delta y} \left[- \text{F}(-1,1) (\text{Phi}(0,1)-\text{Phi}(-1,0))-\text{F}(-1,-1)
			(\text{Phi}(-1,0)-\text{Phi}(0,-1)) \right. \\
			& \qquad \hphantom{\frac{1}{4 \Delta x \Delta y}} \left. + \text{F}(1,-1) (\text{Phi}(1,0)-\text{Phi}(0,-1))+\text{F}(1,1) (\text{Phi}(0,1)-\text{Phi}(1,0)) \right],
		\end{aligned} \\
		& \text{J1} = \frac13 \left(\text{J1pp} + \text{J1px} + \text{J1xp}\right),
	\end{align}
\end{subequations}
we can look at the series expansion of $J_1$:
\begin{subequations}
	\begin{align}
		\mathrm{J1} = & f^{(1,0)}  \phi ^{(0,1)} - f^{(0,1)}  \phi ^{(1,0)} \\
	    & + \frac{1}{6}  \Delta x^2 \left(  f^{(2,0)} \phi^{(1,1)} -  f^{(1,1)} \phi^{(2,0)} -  f^{(2,1)} \phi^{(1,0)} +  f^{(1,0)} \phi ^{(2,1)} + f^{(3,0)} \phi^{(0,1)} -  f^{(0,1)} \phi^{(3,0)} \right) \\
		&+ \frac{1}{6} { \Delta y}^2  \left(f^{(1,0)} \phi^{(0,3)} - f^{(0,3)} \phi^{(1,0)} + f^{(1,1)} \phi^{(0,2)} - f^{(0,2)} \phi^{(1,1)} + f^{(1,2)} \phi^{(0,1)} - f^{(0,1)} \phi^{(1,2)} \right)  \\
	    & + \mathcal{O}(\Delta x^k \Delta y^l \ | \ k+l > 2)
	\end{align}
\end{subequations}
which shows us that the approximation is of order $2$.

\subsubsection*{Fourth Order Scheme (Thirteen-Point-Stencils)}

Similar to above, we define the stencils:
\begin{subequations}
	\begin{align}
		& \begin{aligned}
			\text{J2xp} & = \frac{1}{8 \Delta x \Delta y} \left[ - \text{F}(-1,1)(\text{Phi}(0,2)-\text{Phi}(-2,0)) - \text{F}(-1,-1)(\text{Phi}(-2,0)-\text{Phi}(0,-2)) \right. \\
			& \quad \hphantom{= \frac{1}{8 \Delta x \Delta y}} \left. + \text{F}(1,-1) (\text{Phi}(2,0)-\text{Phi}(0,-2)) + \text{F}(1,1)(\text{Phi}(0,2)-\text{Phi}(2,0)) \right] \, ,
		\end{aligned} \\
		& \begin{aligned}
			\text{J2px} & = \frac{1}{8 \Delta x \Delta y} \left[ - \text{F}(-2,0) (\text{Phi}(-1,1)-\text{Phi}(-1,-1)) + \text{F}(0,-2)(\text{Phi}(1,-1)-\text{Phi}(-1,-1)) \right. \\
			& \quad \hphantom{= \frac{1}{8 \Delta x \Delta y}} \left. - \text{F}(0,2) (\text{Phi}(1,1)-\text{Phi}(-1,1)) + \text{F}(2,0)(\text{Phi}(1,1)-\text{Phi}(1,-1)) \right] \, ,
	    \end{aligned} \\
		& \begin{aligned}
			\text{J2xx} & = \frac{1}{8 \Delta x \Delta y} \left[ (\text{F}(1,1)-\text{F}(-1,-1))(\text{Phi}(-1,1)-\text{Phi}(1,-1)) \right. \\
			& \quad \hphantom{= \frac{1}{8 \Delta x \Delta y}} \left. - (\text{F}(-1,1)-\text{F}(1,-1))(\text{Phi}(1,1)-\text{Phi}(-1,-1)) \right] \, , 
		\end{aligned} \\
		& \text{J2} = \frac{1}{3} (\text{J2px} + \text{J2xp} + \text{J2xx}) \, , \\
        & \text{Jh} = 2\text{J1} - \text{J2} \, ,
	\end{align}
\end{subequations}
and look at the series expansion of $J_h$:
\begin{align}
    \mathrm{Jh} & = f^{(1,0)}(x,y) \phi^{(0,1)}(x,y)-f^{(0,1)}(x,y) \phi ^{(1,0)}(x,y) \\
    &\quad -\frac{1}{12} \Delta x^4 \left( \phi^{(2,1)}(x,y)  f^{(3,0)}(x,y)- f^{(2,1)}(x,y)6\phi^{(3,0)}(x,y)\right)\\
    &\quad +\frac{1}{12} \Delta x^2 \Delta y^2 \left(\phi^{(1,2)}(x,y)  f^{(2,1)}(x,y)-f^{(1,2)}(x,y)  \phi^{(2,1)}(x,y)\right.\\
    &\quad -\frac{1}{36} \Delta x^2 \Delta y^2   \left(  f^{(3,0)}(x,y) \phi^{(0,3)}(x,y)+ f^{(0,3)}(x,y) \phi^{(3,0)}(x,y)\right) \\
    &\quad -\frac{1}{12} \Delta y^4 \left(f^{(1,2)}(x,y) \phi^{(0,3)}(x,y)-f^{(0,3)}(x,y) \phi^{(1,2)}(x,y)\right)\\
    &\quad + \mathcal{O}(\Delta x^k \Delta y^l \ | \ k+l > 4) \, ,
\end{align}
which shows the approximation order of $4$.

\end{document}